\documentclass[pre,dvipsnames,superscriptaddress,onecolumn]{revtex4-2}

\usepackage{amsmath,natbib,amssymb,amsfonts,amsthm,bbm,wrapfig,bm,dsfont,xcolor,graphicx,color}
\usepackage[mathscr]{eucal}
\usepackage[margin=3cm]{geometry}
\newtheorem{theorem}{Theorem}
\newtheorem{corollary}{Corollary}[theorem]
\newtheorem{conjecture}{Conjecture}
\usepackage[colorlinks=true, allcolors=blue]{hyperref}
\allowdisplaybreaks

\begin{document}

\title{The rank and layer distributions in random recursive trees}

\author{Huck Stepanyants}
\affiliation{Department of Physics, Northeastern University, Boston, Massachusetts 02115, USA}
\affiliation{Network Science Institute, Northeastern University, Boston, Massachusetts 02115, USA}
\affiliation{Department of Physics, Harvard University, Cambridge, Massachusetts 02138, USA}
\author{P. L. Krapivsky}
\affiliation{Department of Physics, Boston University, Boston, MA 02215, USA}
\affiliation{Santa Fe Institute, Santa Fe, New Mexico 87501, USA}
\author{Harrison Hartle}
\affiliation{Network Science Institute, Northeastern University, Boston, Massachusetts 02115, USA}
\affiliation{Santa Fe Institute, Santa Fe, New Mexico 87501, USA}
\author{Dmitri Krioukov}
\affiliation{Department of Physics, Northeastern University, Boston, Massachusetts 02115, USA}
\affiliation{Network Science Institute, Northeastern University, Boston, Massachusetts 02115, USA}
\affiliation{Department of Mathematics, Northeastern University, Boston, Massachusetts 02115, USA}
\affiliation{Department of Electrical \& Computer Engineering, Northeastern University, Boston, Massachusetts 02115, USA}

\begin{abstract}
The distribution of node depths in a network is crucial for analyzing network structure. Two measures, rank and layer, quantify how deep inside a network a node is. The rank is the node's distance to the closest leaf, a node of degree~1. The layer is the number of times all leaves must be stripped off for the node to become a leaf. We derive exact recursive expressions for the rank and layer distributions in random recursive trees. We show that the rank and layer distributions decay factorially and geometrically, respectively, and prove the self-averaging of the numbers of nodes with fixed rank or layer. Unlike previous studies, our approach does not depend on labels or a root node, providing a more versatile framework for analyzing rank and layer distributions in complex networks.
\end{abstract}

\maketitle

\section{Introduction}
\label{sec:intro}

Many networks across natural, social, and technological domains have been shown to exhibit a strong core-periphery structure \cite{Holme_2005,Rombach2014}---a densely interconnected central core of nodes, accompanied by an outer periphery where nodes are predominantly connected to the core. The dense core can incubate and amplify perturbations, such as the spread of information or epidemics \cite{everett2000,Sergey06, Dufresne2016,Polanco23,Barber_2015}, which can explosively spill out into the outer peripheral regions. The position of a node in the core-periphery hierarchy is thus a useful predictor of the node's importance and functional role in such dynamical processes. Motivated by this importance, a variety of mechanistic models have been proposed to explain the emergence of such structural organization \cite{borgatti2000models,Verma_2016,Olaizola2024,carmi2007,Rossa_2013}. Several measures of core-periphery structure have been introduced, alongside algorithms for their computation and inference from data \cite{kojaku2017,zhang2015,batagelj2003m,gallagher2021}. For instance, a widely considered quantifier of core-periphery structure in networks has been the $k$-core decomposition~\cite{Sergey06, Ignacio_Alvarez_Hamelin_2008} and its generalizations \cite{doi:10.1137/19M1290607,elliot2020}. Unfortunately, core-periphery metrics are highly nontrivial functions of the global network structure, and hence the progress in the rigorous analysis of these metrics has been relatively slow-paced even in simple network models \cite{fernholz2004cores}. Indeed, even the most basic core-periphery metrics such as the $k$-core decomposition have been analyzed rigorously only for a few network models~\cite{Luczak1991,Pittel1996}.

\begin{figure}
    \centering
    \includegraphics[scale=0.66]{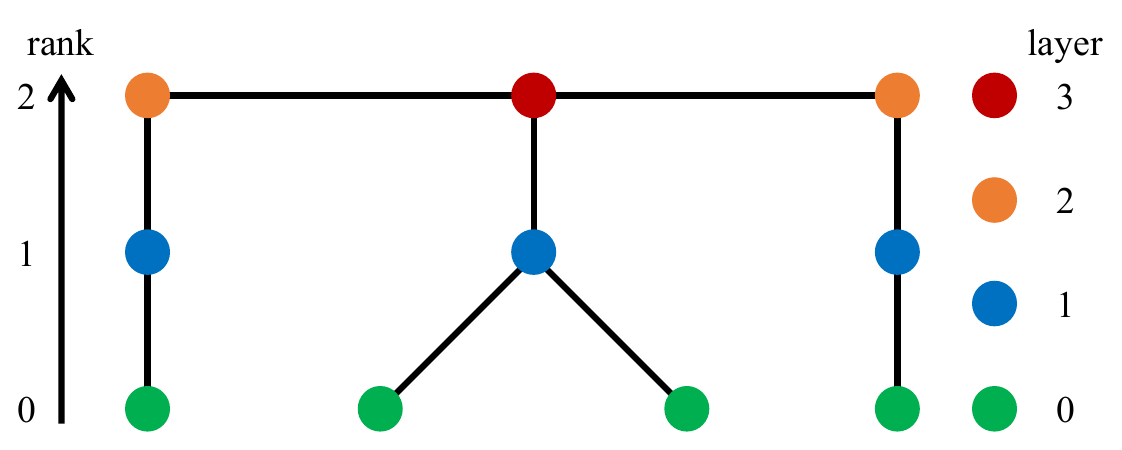}
    \caption{Stratification of a tree with~${n=10}$ nodes by rank and layer. The nodes are arranged vertically in the order of increasing rank (distance to nearest leaf) from~$0$ to~$2$. Node layer (number of rounds of leaf-removal to become a leaf) is indicated by color, from green, layer~0, to red, layer~3. In any tree, the rank and layer of leaves are both equal to~$0$, while the rank of their nonleaf neighbors is~$1$.}
    \label{fig:rank_and_layer}
\end{figure}

Here, we focus on two measures of core-periphery structure---the \emph{rank} and the \emph{layer}---in the context of the simplest growing network model, the random recursive tree (RRT). The RRT model was introduced long ago \cite{RRT-Tapia,RRT-Moon,RRT-Pittel} and it acquired the status of a null model for growing networks \cite{Drmota,Newman-book,KRB,Frieze}. In the RRT, nodes arrive one at a time, each connecting to one existing node selected uniformly at random. The rank and layer both stratify nodes in a network by their distance from the leaves: the rank of a node is the distance from it to the nearest leaf; the layer of a node is the number of steps needed to make it a leaf by recursive rounds of leaf removals; see Fig.~\ref{fig:rank_and_layer}. The layer decomposition is related to the ${k}$-core decomposition \cite{SEIDMAN1983269,Pittel1996,Sergey06}, even though the $k$-core decomposition of any tree is trivial: all nodes are in the $1$-core. The layer decomposition is also the special case of the onion decomposition~\cite{Dufresne2016}, which generalizes the $k$-core decomposition---namely, the layer is the onion depth into the $1$-core. We are primarily interested in the rank and layer distributions, i.e., the probabilities that a node chosen uniformly at random in the network has rank~$k$ or is at layer~$k$, which we denote by~${r_k}$ and~${\ell_k}$, respectively.

Thanks to the simplicity of RRTs, they have received rigorous treatments for a variety of sophisticated metrics, such as tree height and width \cite{Drmota,Flajolet,Frieze,Devroye87,RRT-Pittel}. However, the rank and layer distributions in RRTs have not been previously obtained. The ${\ell}$-protectedness, a property closely related to rank, has been studied extensively in several random trees \cite{Mahmoud12,Mahmoud15,Copenhaver17}. A node is~${\ell}$-protected if the shortest path from it to any {\it descendant} leaf is of length~${\ell}$. The protectedness of a node is thus a function of the subtree that it roots. Using a result by Aldous~\cite{Aldous1991} on arbitrary properties of subtrees, Devroye and Janson~\cite{Devroye14} computed the asymptotic fraction of ${\ell}$-protected nodes in RRTs, Galton-Watson trees, and binary search trees. Additionally, these fractions have been computed for the first few values of~${\ell}$ in 1-2 trees~\cite{Bona19}, random increasing plane trees~\cite{Bona21}, and digital search trees~\cite{Prodinger12}. However, since a node's  rank and layer depend not only on its subtree but on the entire tree, the methods used to study $\ell$-protectedness are not immediately applicable here.

Our main results are stated in Sec.~\ref{sec:results}. The key techniques and ideas behind the proofs---tree-splitting, node anchoring, and a recursive method---are presented in Sec.~\ref{sec:recursive}. The proofs themselves are documented in Sec.~\ref{sec:proofs}. Appendix \ref{ssec:large_k_proof_outline} outlines the traveling wave method used to arrive at our conjecture for the tail of the layer distribution. Here is a brief summary of our main results:
\begin{itemize}
  \item Exact recursive expressions for the rank and layer distributions $r_k$ and $\ell_k$ in infinite RRTs (Theorems~\ref{the:ak} and \ref{the:pk}) and their evaluations for small values of~$k$ (Corollary~\ref{cor:smallk} and Table~\ref{tab:table}).
  \item Characterization of large-$k$ asymptotics of $r_k$ and $\ell_k$ (Theorems~\ref{the:rank_dist_tail} and~\ref{the:layer_dist_tail} and Conjecture~\ref{con:layer_dist_tail}).
  \item Proof of self-averaging of the numbers of nodes with fixed rank or layer, and a lower bound for the speed of convergence of the expected values of these fractions to their infinite-size limits $r_k$ and $\ell_k$ (Theorems~\ref{the:Rnk_variance} and~\ref{the:Rnk_expectation}).
  \item Exact expressions for $r_0(n)=\ell_0(n)$, $r_1(n)$, and $\ell_1(n)$ in RRTs of any finite size  $n$ (Theorem~\ref{the:finite_n}).
\end{itemize}

\section{Main results}
\label{sec:results}

Before stating our results, we introduce some \textbf{notations and definitions}. A random recursive tree (RRT) is a labeled random tree, with nodes labeled by their arrival times $\{1,2,\ldots,n\}=[n]$. Upon its arrival, a new node connects to one existing node chosen uniformly at random from the set of existing nodes. In an RRT $G_n$ of size $n$, let~$R_{n,k}$ denote the number of nodes with rank~$k$, and let~$L_{n,k}$ denote the number of nodes at layer~$k$, while $r_k(n)=\langle R_{n,k}/n \rangle_{G_n}$ and $\ell_k(n)=\langle L_{n,k}/n \rangle_{G_n}$ denoting the expected fractions of nodes with rank~$k$ or at layer~$k$, respectively. Let~${r_k=\lim_{n\rightarrow\infty} r_k(n)}$ and~${\ell_k=\lim_{n\rightarrow\infty} \ell_k(n)}$ denote the $n\to\infty$ limits of these fractions (they do exist).

Define $\bar{r}_k(n)=\sum_{k'=k}^\infty r_{k'}(n)$ and $\bar{r}_k=\sum_{k'=k}^\infty r_{k'}$ to be the cumulative distributions of expected rank fractions in finite-size and infinite-size RRTs, respectively. Similarly, define $\bar{\ell}_k(n)=\sum_{k'=k}^\infty \ell_{k'}(n)$ and $\bar{\ell}_k=\sum_{k'=k}^\infty \ell_{k'}$ to be the cumulative distributions of the expected layer fractions in finite-size and infinite-size RRTs, respectively. Trivially, $\bar{r}_0(n)=\bar{r}_0=\bar{\ell}_0(n)=\bar{\ell}_0=1$. Finally, the relation ``$\lesssim$'' in $f(n)\lesssim g(n)$ means that there exists~${C}$ such that~${f(n)/g(n) \leq C}$ for all~${n}$, and all $\log$s are natural.\\

Our \textbf{main results}, Theorems~\ref{the:ak} and \ref{the:pk}, are about the cumulative rank and layer distributions in infinite RRTs.

\begin{theorem}
\label{the:ak}
The cumulative rank distribution is given by 
\begin{equation}
    \label{eq:ak1}
    \bar{r}_k = \int\displaylimits_{0 < \alpha_1 < \dots < \alpha_k < 1} \alpha_1 f_k(1-\alpha_k) \prod\displaylimits_{i=1}^{k-1} \left[1 + f_i(1-\alpha_i) \right]\, d \alpha_1 \ldots d \alpha_k 
\end{equation}
with convention $\prod_{i=1}^{k-1}=1$ when $k=1$ and functions $f_i(x)$ defined recursively via 
\begin{equation}
    \label{eq:ak2}
    f_0(x) = \frac{1}{1-x},
    \quad f_i(x) = \exp\!\left[\int_0^x f_{i-1}(x^\prime) \, dx^\prime \right] - 1 \quad \forall i \geq 1.
\end{equation}
\end{theorem}
Theorem~\ref{the:ak} recovers previously obtained values~${\bar{r}_1 = 1/2}$ and~${\bar{r}_2 = 1/2 - 1/e}$~\cite{Devroye14,Mahmoud15}.

\begin{figure}
    \centering
    \includegraphics[width=0.8\textwidth,trim=0 30 0 0,clip]{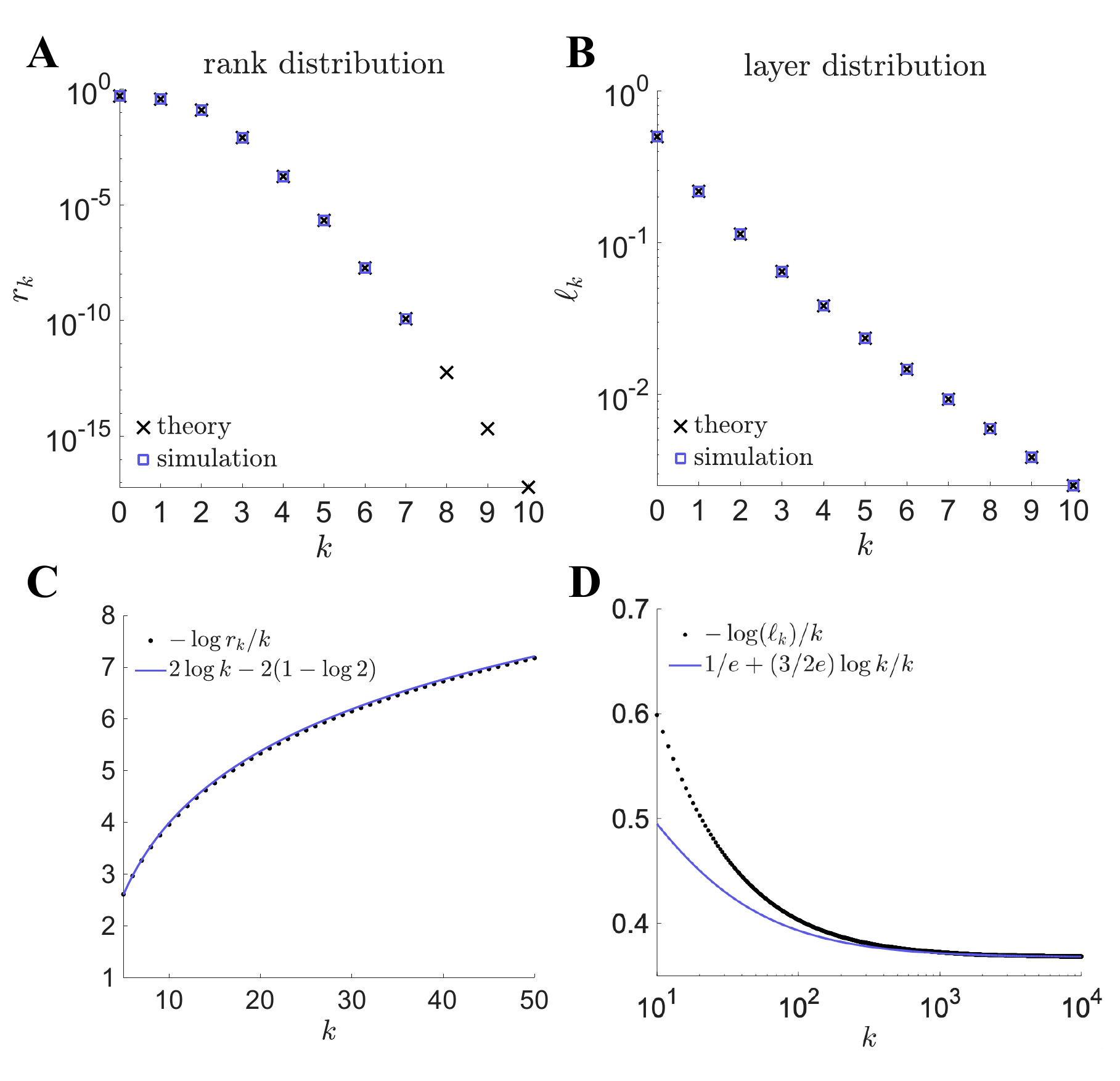}
    \caption{\textbf{(A)}~The numeric evaluation of the rank distribution values~$r_k=\bar{r}_k-\bar{r}_{k+1}$ from Theorem~\ref{the:ak} versus the corresponding empirical values averaged across~${1.5 \times 10^9}$ simulated RRTs of size~${n=10,000}$.  \textbf{(B)}~The numeric evaluation of the layer distribution values~$\ell_k=\bar{\ell}_k-\bar{\ell}_{k+1}$ from Theorem~\ref{the:pk} versus the corresponding empirical values for the same simulated RRTs as in~\textbf{(A)}.  \textbf{(C)}~The numeric evaluation of $-\log(r_k)/k$ from Theorem~\ref{the:ak} (the black dots) versus the large-${k}$ asymptotic of $r_k$ from Theorem~\ref{the:rank_dist_tail} (the blue curve). \textbf{(D)}~The numeric evaluation of $-\log(\ell_k)/k$ from Theorem~\ref{the:pk} (the black dots) versus the large-${k}$ asymptotic of $\ell_k$ from Conjecture~\ref{con:layer_dist_tail} (the blue curve).}
    \label{fig:RLdist_numerical}
\end{figure}

\begin{theorem}
\label{the:pk}
    The cumulative layer distribution is given by
\begin{equation}
    \label{eq:pk1}
    \bar{\ell}_k = \displaystyle\int_0^\infty  s_k(y) \, e^{-y} dy,
\end{equation}
with functions $s_i(y)$ defined recursively via 
\begin{equation}
    \label{eq:pk2}
    s_0(y) = 1,
    \quad s_i(y) = 1 - \exp\!\left[ -\int_0^y s_{i-1}(y^\prime) \, dy^\prime \right] \quad \forall i \geq 1.
\end{equation}
\end{theorem}
Theorem~\ref{the:pk} implies that~${\bar{\ell}_1 = 1/2}$ (reflecting the asymptotic leaf fraction), as well as~${\bar{\ell}_2 = 3-e}$. The proofs of Theorems~\ref{the:ak} and~\ref{the:pk} are in Secs.~\ref{ssec:recursive_method}, \ref{ssec:ak_proof}, and~\ref{ssec:pk_proof}, and the simulation results are in Fig.~\ref{fig:RLdist_numerical}(A,B). The proofs of Theorems~\ref{the:ak} and~\ref{the:pk} not merely give the cumulative rank and layer distributions in infinite RRTs, but also show the convergence, $\bar{r}_k(n)\to \bar{r}_k$ and $\bar{\ell}_k(n) \to \bar{\ell}_k$, in the $n\to\infty$ limit for any fixed $k$. The upper bounds for the rate of convergence are given in Theorem~\ref{the:Rnk_expectation}. 

Although the closed-form solutions of Eqs.~(\ref{eq:ak1}-\ref{eq:pk2}) are not obtainable for all~$k$, we evaluate them for small values of~$k$, resulting in the following corollary.
\begin{corollary}
\label{cor:smallk}
For $k=0,1,\ldots,8$, the expected values of the rank $r_k=\bar{r}_k-\bar{r}_{k+1}$ and layer~$\ell_k=\bar{\ell}_k-\bar{\ell}_{k+1}$ fractions are as shown in Table~\ref{tab:table}.
\end{corollary}
\begin{table}
    \begin{tabular}{ |c|c|c|c|c|c|c|c|c|c|c|c| }
     \hline
     ${\pmb{k}}$ & 0 & 1 & 2 & 3 & 4 & 5 & 6 & 7 & 8 \\
     \hline
     ${\pmb{r_k}}$ & ${1/2}$ & ${1/e}$ & ${0.124}$ & ${0.00805}$ & ${1.68 \times 10^{-4}}$ & ${2.14 \times 10^{-6}}$ & ${1.88 \times 10^{-8}}$ & ${1.20 \times 10^{-10}}$ & ${5.76 \times 10^{-13}}$ \\
     \hline
     ${\pmb{\ell_k}}$ & ${1/2}$ & ${e-5/2}$ & ${0.114}$ & ${0.0647}$ & ${0.0384}$ & ${0.0235}$ & ${0.0146}$ & ${0.00928}$ & ${0.00595}$ \\
     \hline
    \end{tabular}
    \caption{Values of~${r_k}$ and~${\ell_k}$ for~${k \leq 8}$ computed from Eqs.~(\ref{eq:ak1}-\ref{eq:pk2}).}
\label{tab:table}
\end{table}
The numbers in Table~\ref{tab:table} are calculated exactly for~${k=0,1,2}$ in Sec.~\ref{ssec:rl_small_k}, while for $k=3,\ldots,8$, they are obtained by solving the recursions in Theorems~\ref{the:ak} and~\ref{the:pk} numerically.

The other way around, for large values of~$k$, we show that the tail of the rank distribution~$r_k$ decays factorially in Theorem~\ref{the:rank_dist_tail}, and that the tail of the layer distribution~$\ell_k$ is geometric with an exponent that we can compute in Theorem~\ref{the:layer_dist_tail}.
\begin{theorem}
\label{the:rank_dist_tail}
For $k\gg 1$, the rank distribution decays as
\begin{equation}
\label{eq:rank_tail}
    r_k = \exp\!\left[ -2k \log k + 2(1-\log 2)k + O(\log k) \right].
\end{equation}
\end{theorem}
The proof is in Sec.~\ref{ssec:rank_tail_proof}, and the numeric evidence is in Fig.~\ref{fig:RLdist_numerical}(C).

\begin{theorem}
\label{the:layer_dist_tail}
The cumulative layer distribution~${\bar{\ell}_k}$ is asymptotically bounded for large~${k}$ by exponential functions with the decay rates~${1}$ and~${1/e}$,
\begin{equation}
\label{eq:layer_bounds}
    e^{-k} \lesssim \bar{\ell}_k \lesssim e^{-k/e}.
\end{equation}
\end{theorem}
The proof is in Sec.~\ref{ssec:layer_tail_proof}.

We further conjecture the stronger statement that the layer distribution~$\ell_k$ attains the upper bound in Eq.~\eqref{eq:layer_bounds}.
\begin{conjecture}
\label{con:layer_dist_tail}
The layer distribution~$\ell_k$ is asymptotically
\begin{equation}
\label{eq:layer_tail}
    \ell_k = \exp\!\left[ -\frac{1}{e} k - \frac{3}{2e} \log k - O(1) \right].
\end{equation}
\end{conjecture}
In Appendix~\ref{ssec:large_k_proof_outline} we provide analytic and numerical evidence to support this conjecture, the latter is also presented in Fig.~\ref{fig:RLdist_numerical}(D).

Theorems~\ref{the:rank_dist_tail} and~\ref{the:layer_dist_tail} suggest asymptotics for the maximal rank~${R}$ and maximal layer~${L}$ across all nodes in an RRT of size~${n \gg 1}$. Indeed, by setting~${r_k n \sim 1}$ and~${\ell_k n \sim 1}$, one obtains~${2R\log R \approx \log n}$ and~${\frac{1}{e} L + \frac{3}{2e} \log L \approx \log n}$, respectively. Solving these equations leads to the following 
\begin{conjecture}
\label{con:max}
The maximum rank~$R$ and layer~$L$ in an RRT of size~$n$ are asymptotically
\begin{align}
  R &\approx \frac{1}{2} \frac{\log n}{\log \log n}, \\
  L &\approx e\log n - \frac{3}{2} \log \log n.
\end{align}
\end{conjecture}

\begin{figure}
    \centering
    \includegraphics[width=0.8\textwidth]{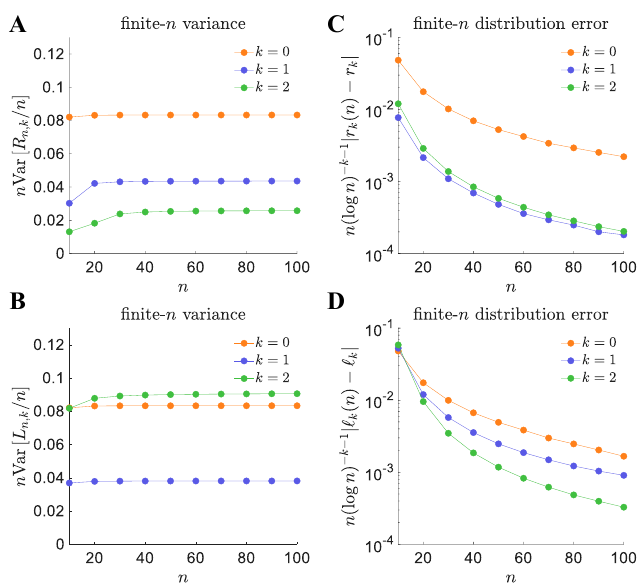}
    \caption{\textbf{(A, B)} The variances of the rank and layer fractions rescaled by~${n}$, ${n\text{Var}\left[R_{n,k}/n\right]}$ and~${n\text{Var}\left[L_{n,k}/n\right]}$ (dots connected by lines), in ${2\times10^8}$ simulated RRTs for each value of~$n$ and $k=0,1,2$, are upper bounded by constants in agreement with Theorem~\ref{the:Rnk_variance}.
    \textbf{(C, D)}~The differences between the expected rank and layer fractions~${r_k(n)}$ and~${\ell_k(n)}$ and their corresponding limits rescaled by~${n (\log n)^{-k-1}}$ (dots connected by lines) in the same RRTs as in \textbf{(A, B)} are upper bounded by constants in agreement with Theorem~\ref{the:Rnk_expectation}.}
    \label{fig:variance_convergence_numerical}
\end{figure}

The knowledge of limiting expected values of particular properties---$r_k$ and $\ell_k$ in our case---may be next to useless for characterizing individual RRT realizations unless we also know that the distributions are sharply peaked at their expected values. The next two Theorems~\ref{the:Rnk_variance} and \ref{the:Rnk_expectation} establish two concentration results---self-averaging and the speed of convergence---for the fractions of nodes of rank~$k$, $R_{n,k}/n$, and at layer~$k$, $L_{n,k}/n$.
\begin{theorem}
\label{the:Rnk_variance}
The fractions of nodes of rank ${k}$ and at layer $k$ are asymptotically self-averaging, meaning that the following asymptotic inequalities hold for any fixed~$k$:
\begin{equation}
\begin{aligned}
\label{eq:self_avg}
    \mathrm{Var} & \left[ \frac{R_{n,k}}{n} \right] \lesssim \frac{1}{n}, \\
    \mathrm{Var} & \left[ \frac{L_{n,k}}{n} \right] \lesssim \frac{1}{n}.
\end{aligned}
\end{equation}
\end{theorem}
The proof is in Sec.~\ref{ssec:variance_proof}, and the numeric evidence is in Fig.~\ref{fig:variance_convergence_numerical}(A,B). 

\begin{theorem}
\label{the:Rnk_expectation}
The differences between the expected fractions of nodes of rank~$k$ or at layer~${k}$ in RRTs of size~$n$, and their respective values in the infinite RRTs, obey the following asymptotic inequalities for any fixed~$k$:
\begin{equation}\begin{aligned}
\label{eq:self_avg_2}
    \left| r_k(n) - r_k \right|
    &\lesssim \frac{(\log n)^{k+1}}{n}, \\
    \left| \ell_k(n) - \ell_k \right| &\lesssim \frac{(\log n)^{k+1}}{n}.
\end{aligned}\end{equation}
\end{theorem}
The proof is in Sec.~\ref{ssec:Rnk_expectation}, and the numeric evidence is in Fig.~\ref{fig:variance_convergence_numerical}(C,D). 

\begin{figure}
    \centering
    \includegraphics[width=1\textwidth]{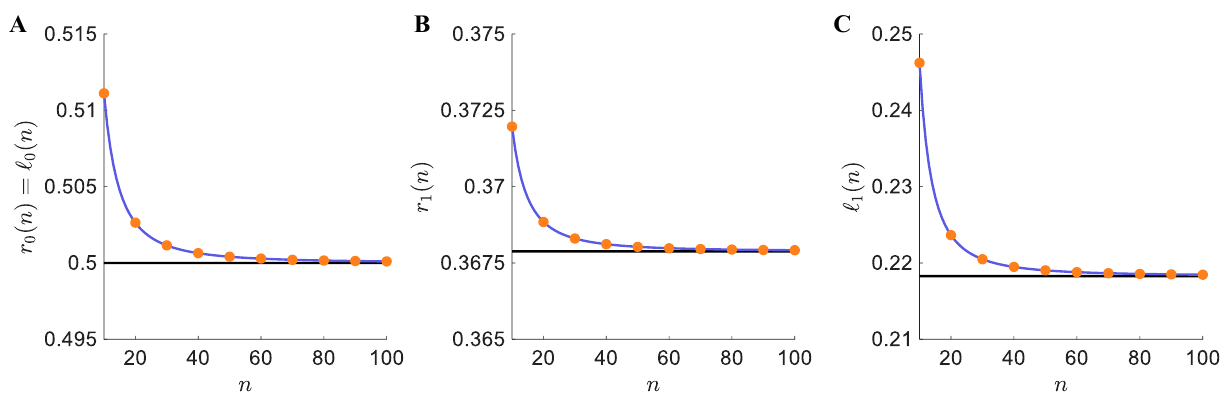}
    \caption{The average fractions of leaves~\textbf{(A)}, nodes of rank~1~\textbf{(B)}, and at layer~1~\textbf{(C)} across~${2\times10^8}$ simulated RRTs for each value of~${n}$ (the orange dots) versus the corresponding expected fractions in Theorem~\ref{the:finite_n} (the blue curves) and their $n\to\infty$ limits in Corollary~\ref{cor:smallk} (the black lines). For each simulation point, the uncertainty is smaller than the marker size.}
    \label{fig:smalln_numerical}
\end{figure}

Finally, the explicit formulae for the expected fractions of nodes of a given rank or at layer~${k=0,1}$ in an RRT of any finite size~$n$ are:
\begin{theorem}
\label{the:finite_n}
For~${n>2}$, the expected fraction of leaves in an RRT of size~${n}$ is
\begin{equation}
\label{r-0n}
    r_0(n) = \ell_0(n) = \frac{1}{2} + \frac{1}{n(n-1)}\,.
\end{equation}
The expected fraction of nodes of rank~${1}$ and at layer~$1$ are, respectively,
\begin{align}
\label{r-1n}
    r_1(n) =& \sum_{j=0}^{n-2} \frac{(-1)^j}{j!}  + \frac{1}{n(n-1)} \sum_{j=0}^{n-3} \frac{(-1)^j}{j!} (j+2)\\
\label{r-1n:asymp}
    =&~e^{-1} + \frac{e^{-1}}{n(n-1)} + O\!\left(\frac{1}{(n+1)!}\right),\\
\label{ell-1n}
    \ell_1(n) =& -\frac{1}{2} + \frac{1}{n} - \frac{1}{n-1}
    + \frac{2}{n!} + \sum_{v=2}^{n-1} \frac{1}{v!} - \frac{1}{n(n-1)} \sum_{v=0}^{n-3} \frac{1}{v!} + \frac{2}{n!} \sum_{j=1}^{n-2} \sum_{v=0}^{j-1} \frac{j!}{v!} \\
\label{l-1n:asymp}
    =&~e - \frac{5}{2} + \frac{e-1}{n^2} + O \left( \frac{1}{n^3} \right)\,.
\end{align}
\end{theorem}
The proof is in Sec.~\ref{ssec:finite_n}, and the numeric evidence is in Fig.~\ref{fig:smalln_numerical}.

The fractions $r_1(n)$ and ${\ell_1(n)}$ in Eqs.~\eqref{r-1n} and~\eqref{ell-1n} are positive for all $n>2$. The explicit results in Eqs.~\eqref{r-0n}, \eqref{r-1n:asymp}, and \eqref{l-1n:asymp} show that~$r_0(n)$, ${\ell_0(n)}$, ${r_1(n)}$, and~${\ell_1(n)}$ converge to their limits with rates~$O(n^{-2})$, that is, faster than the general upper bounds, $n^{-1}(\log n)^{k+1}$ in Theorem~\ref{the:Rnk_expectation}, so strengthening these upper bounds may be feasible. The number of leaves in an RRT, $R_{n,0}$ or equivalently $L_{n,0}$, is a simple random quantity whose mean value~\eqref{r-0n} as well as higher moments and cumulants are known~\cite{hartle26statistics}. More detailed characterizations of the random quantities $R_{n,1}$ and $L_{n,1}$ are interesting open problems.

\section{Main techniques and ideas behind the proofs}
\label{sec:recursive}

In this section, we sketch the main ideas behind the proofs of our main results. The key techniques are tree splitting and anchoring that we describe in Sec.~\ref{ssec:tree_splitting_and_anchoring}. We use them in our recursive method, to illustrate which we outline the proof of Theorem~\ref{the:ak} in Sec.~\ref{ssec:recursive_method}. The proof of Theorem~\ref{the:pk} is very similar.

\subsection{Tree splitting and anchoring}
\label{ssec:tree_splitting_and_anchoring}

\begin{figure}
    \centering
    \includegraphics[width=0.66\textwidth]{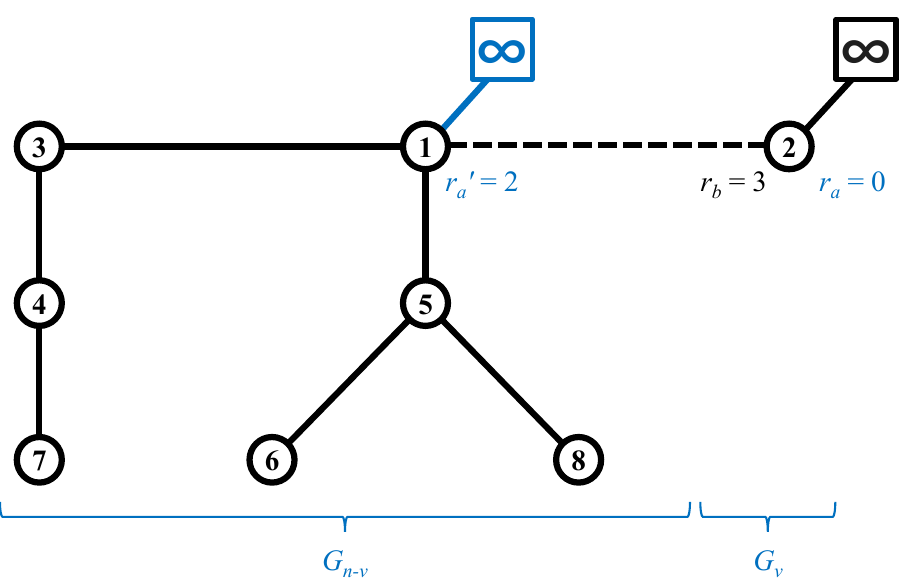}
    \caption{Example of tree-splitting and anchoring in an RRT with~${n=8}$ nodes. The anchor shown by the black boxed~$\infty$ is attached to node ${t=2}$. The rank of~$t$ becomes $r_b=3$. Then the edge between ${t=2}$ and its ancestor ${t^\prime=1}$ is deleted, after which ${t^\prime}$ is anchored as well. The rank of~$t$ is now $r_a=0$, while the rank of $t^\prime$ is now $r_a^\prime=2$. The elements of the original tree are colored black, while changes resulting from this tree splitting process are colored blue. In the figure, the size of the split subtree $G_v$ rooted at~${t=2}$ is ${v=1}$, corresponding to Case~1 in Sec.~\ref{ssec:recursive_method}.}
    \label{fig:splitting_anchoring}
\end{figure}

One key technique we employ in the proofs herein is tree splitting. This technique was used in Ref.~\cite{Mahmoud15}. Tree splitting is the deletion of the edge connecting a node with label~${t \geq 2}$ to its ancestor in an RRT~$G_n$ of size~${n}$. This deletion creates two smaller trees, the subtree~$G_v$ of size~$v$ rooted at node~${t}$, and the remaining part~$G_{n-v}$ of~$G_n$, as illustrated in Fig.~\ref{fig:splitting_anchoring}. Both $G_v$ and $G_{n-v}$, conditioned on their sizes, are also RRTs.

Consider a node ${s \neq \{1,t\}}$. This node is added at time ${s}$ and is attached to one of the ${s-1}$ existing nodes in the RRT chosen uniformly at random. If~${s}$ is attached to~$G_v$, then the label of the node to which it is attached is distributed uniformly amongst the labels of all the nodes in~$G_v$. Thus, the growth mechanism of~$G_v$ is the same as that of an RRT. Similarly, if ${s}$ is attached to~$G_{n-v}$, then the label of the node to which it is attached is distributed uniformly amongst the labels of all the nodes in~$G_{n-v}$. Therefore, both $G_v$ and $G_{n-v}$ are RRTs. However, we should keep in mind that the node label sets of $G_v$ and $G_{n-v}$ are two nonoverlapping subsets of $[n]$,  rather than $[v]$ and $[n-v]$ if the labeling was standard.

As we show in Sec.~\ref{ssec:subtree_size_distribution}, the distribution of the size~$v\in[n-t+1]$ of~$G_v$ is
\begin{equation}
\label{eq:sntm}
        h_{n,t}(v) =\frac{\binom{n-v-1}{t-2}}{\binom{n-1}{t-1}}.
\end{equation}
Of particular interest is~$t=2$, in which case this distribution is uniform:
\begin{equation}
\label{unif_sn2}
    h_{n,2}(v)=\frac{1}{n-1}.
\end{equation}

Another key technique we rely on to derive our recursions for the rank and layer is \emph{anchoring}. Anchors are infinitely long chains of nodes, shown by the boxed~$\infty$s in Fig.~\ref{fig:splitting_anchoring}. Anchors are not parts of an RRT, but by attaching them to an RRT, we can control the ranks of nodes in it. In particular, by attaching anchors to the two nodes directly connected by an edge that was deleted during a tree splitting event, we ensure that no new leaves are created during the event. Note that attaching an anchor to a nonleaf in an RRT will have no effect on the ranks of any of its nodes.

\subsection{The recursive method}
\label{ssec:recursive_method}

Here we present the most important steps in the proof of Theorem~\ref{the:ak}; the remaining calculations are covered in Secs.~\ref{ssec:ak_proof} and~\ref{ssec:pk_proof}. In our recursive method, the two main players are:
\begin{itemize}
\item[1)] $b_{n,t}^k$: the probability that node~$t$ has rank~$\ge k$;
\item[2)] $c_{n,t}^k$: the same as $b_{n,t}^k$, except that $t$ is anchored.
\end{itemize}

The method proceeds as follows. We first anchor node~${t \geq 2}$ and denote its rank after anchoring by~${r_b}$. In the illustration in Fig.~\ref{fig:splitting_anchoring}, $t=2$ and $r_b=3$. We then split $G_n$ at $t$ by removing $t$'s link to its ancestor in~$G_n$ whose label we denote by~$t^\prime$, which is uniformly distributed on~$[t-1]$, and then we anchor~${t^\prime}$ as well. Node~$t$ is now the root of $G_n$'s subtree~$G_v$ of size~$v$. We denote by~$r_a$ and~$r_a^{\prime}$, respectively, the new ranks of~${t}$ and~$t^\prime$. In Fig.~\ref{fig:splitting_anchoring}, $r_a=0$ and $r_a^\prime=2$.

We now consider the following two cases:
\begin{enumerate}
\item if node~${t}$ was a leaf, then~${v=1}$ and~${r_b = r_a^{\prime} + 1}$, so~${r_b \ge k}$ if and only if~${r_a^{\prime} \ge k-1}$;
\item if node~${t}$ was a nonleaf, then~${2 \leq v \leq n-t+1}$ and~${r_b = \min \{ r_a,r_a^{\prime}+1 \}}$, so ${r_b \ge k}$ if and only if~${r_a \ge k}$ and ${r_a^{\prime} \ge k-1}$.
\end{enumerate}
In view of these observations, we deduce that the tree splitting at node~${t}$ leads to the recursion
\begin{align}
\label{eq:r_outline1}
    & c_{n,t}^k = 
    h_{n,t}(1) ~\mathbbm{E}_{t^\prime} \left[ c_{n-1,t'}^{k-1} \right] + \sum_{v=2}^{n-t+1} h_{n,t}(v) ~\mathbbm{E}_{t'<t} \left[ c_{n-v,t'}^{k-1} \, c_{v,1}^k \right],
\end{align}
with the initial condition ${c_{n,1}^k = c_{n,2}^k}$, which follows from the symmetry of the RRT distribution under permutations of the first two nodes.

Next, we use tree splitting again to derive a formula for the~${b_{n,t}^k}$ in terms of the~${c_{n,t}^k}$. This time, the initial node~${t}$ is initially unanchored. After removing the edge between~${t}$ and its ancestor~${t'}$, both nodes are then anchored, as before. We define~${v}$, ${r_b}$, ${r_a}$, and~${r_a'}$ as previously, and find
\begin{equation*}
r_b = \min\{r_a,r_a'+1\}
\end{equation*}
regardless of~${v}$. Therefore, ${r_b \ge k}$ if and only if~${r_a \ge k}$ and~${r_a' \ge k-1}$. This leads to the formula
\begin{align}
    \label{eq:r_outline5}
    & b_{n,t}^k = \sum_{v=1}^{n-t+1} h_{n,t}(v) ~\mathbbm{E}_{t'<t} \left[ c^{k-1}_{n-v,t'} \, c_{v,1}^k \right].
\end{align}

Equations~\eqref{eq:r_outline1} and~\eqref{eq:r_outline5} completely determine~${b_{n,t}^k}$, and are solved in Sec.~\ref{ssec:ak_proof} for~${n\rightarrow\infty}$ using the method of generating functions. The final step in the proof is to compute~${\bar{r}_k}$, which is just the expected value of~${b_{n,t}^k}$ over~${t = 1,2,\dots,n}$ in the limit~${n \rightarrow \infty}$,
\begin{align}\label{eq:r_outline4}
    & \bar{r}_k = \lim_{n \rightarrow \infty} \frac{1}{n} \sum_{t=1}^n b_{n,t}^k.
\end{align}

The proof of Theorem~\ref{the:pk} for the layer distribution is similar, although the process of tree splitting and anchoring is slightly more involved. First, define~$b_{n,t}^k$, $c_{n,t}^k$, ${\ell_b}$, ${\ell_a}$, and ${\ell_a'}$ analogously to before. That is, ${\ell_b}$ is the layer of the anchored node~$t$ before tree splitting, ${\ell_a}$ is its layer after tree splitting, ${\ell_a'}$ is the layer of its ancestor~$t'$ after tree splitting, while
\begin{itemize}
\item[1)] $b_{n,t}^k$ is the probability that node~$t$ is at layer~$\ge k$,
\item[2)] $c_{n,t}^k$ is the same as $b_{n,t}^k$, except that $t$ is anchored.
\end{itemize}

We now consider the case where~${t}$ is initially anchored. Then~${\ell_b = \max \{ \ell_a,\ell_a'+1 \}}$, so~${\ell_b \ge k}$ if and only if~${\ell_a \ge k}$ or~${\ell_a' \ge k-1}$. This leads to the recursion
\begin{align}\label{eq:l_outline1}
    & c_{n,t}^k = \sum_{v=1}^{n-t+1} h_{n,t}(v) ~\mathbbm{E}_{t'<t} \left[ 1 - (1 - c_{n-v,t'}^{k-1})(1 - c_{v,1}^k) \right].
\end{align}
Now, consider the case where~${t}$ is initially unanchored. We consider the following two subcases. The first one is: (1)~if~${\ell_a' \geq k-1}$, then~${\ell_b \geq k}$ if and only if~${\ell_a \geq k}$. We note that in this subcase, the anchor we attach to~${t}$ takes the place of the ancestral branch from~${t}$ to~${t'}$. The second subcase is: (2)~if~${\ell_a' < k-1}$, then~${\ell_b \geq k}$ if and only if the layer of the \emph{unanchored} node~${t}$ \emph{after} tree splitting is~${\geq k}$. In this case, we ignore the ancestral branch from~${t}$ to~${t'}$ entirely, since it had no effect on whether~${r_a \geq k}$. The resulting recursion is
\begin{align}\label{eq:l_outline2}
    & b_{n,t}^k = \sum_{v=1}^{n-t+1} h_{n,t}(v) ~\mathbbm{E}_{t'<t} \left[ c_{n-v,t'}^{k-1} \, c_{v,1}^k + \left( 1 - c_{n-v,t'}^{k-1} \right) b_{v,1}^k \right].
\end{align}

This recursion simplifies drastically in the~${n \rightarrow \infty}$ limit. Indeed, since the height of an RRT scales logarithmically with the size of the tree \cite{RRT-Pittel}, an RRT of size $n$ contains a branch of length~${\sim \log n}$ with high probability. This implies that the layer of any anchored node in such a tree is~${\sim \log n}$, and therefore~${\lim_{n \rightarrow \infty} c_{n-v,t'}^{k-1} = 1}$ in Eq.~\eqref{eq:l_outline2}, so the recursion becomes
\begin{align}
\label{eq:l_outline3}
    & b_{n,t}^k = \sum_{v=1}^{n-t+1} h_{n,t}(v) c_{v,1}^k.
\end{align}

Equations~\eqref{eq:l_outline1} and~\eqref{eq:l_outline3} completely determine~${b_{n,t}^k}$, and are solved in Sec.~\ref{ssec:pk_proof} for~${n\rightarrow\infty}$. The final step in the proof is to compute~${\bar{\ell}_k}$, which is given by the same formula as for~${\bar{r}_k}$, see Eq.~\eqref{eq:r_outline4}.

Alternatively, without taking any limits, we can solve the recursions~\eqref{eq:r_outline1}, \eqref{eq:r_outline5}, \eqref{eq:l_outline1}, \eqref{eq:l_outline2} exactly in the cases~${(k=2, ~t=1)}$ and~${(k=1, ~t \geq 1)}$ to obtain the exact finite-${n}$ results in Theorem~\ref{the:finite_n}. The details are in Sec.~\ref{ssec:finite_n}.

\section{Proofs}
\label{sec:proofs}
Here we provide the remaining proofs and calculations not covered in the earlier sections.

\subsection{Subtree size distribution}
\label{ssec:subtree_size_distribution}

To determine ${h_{n,t}(v)}$, the probability that the subtree of node~${t}$ has size~${v}$ in an RRT of size~${n}$, we first calculate ${N_{n,t}(v)}$, the number of RRTs of size ${n}$ for which the subtree of node ${t}$ has size ${v}$. There are~${(t-1)!}$ ways to attach the first~${t}$ nodes of the RRT. Then out of the remaining~${n-t}$ nodes, there are~${\binom{n-t}{v-1}}$ ways to choose which ones will make up the subtree and which will make up the remaining part of the tree. The subtree and the remaining part of the tree can be assembled in~${(v-1)!}$ and ${(t-1) \cdot t \cdots (n-v-1)}$ ways, respectively. So we have
\begin{equation}\begin{aligned}
    N_{n,t}(v) & = (t-1)! \binom{n-t}{v-1} (v-1)! \left[ (t-1) \cdot t \cdots (n-v-1) \right] \\
    & = \binom{n-v-1}{t-2} (t-1)! (n-t)!.
\end{aligned}\end{equation}
Then to obtain~${h_{n,t}(v)}$, we normalize~${N_{n,t}(v)}$ by the total number of RRTs of size~${n}$, which is ${(n-1)!}$:
\begin{equation}\begin{aligned}
    h_{n,t}(v) &= \frac{1}{(n-1)!} N_{n,t}(v) \\
    &= \frac{1}{(n-1)!} \binom{n-v-1}{t-2} (t-1)! (n-t)! \\
    &= \binom{n-1}{t-1}^{-1} \binom{n-v-1}{t-2}.
\end{aligned}\end{equation}

In the scaling limit
\begin{equation}
\label{scaling}
n \rightarrow \infty, \quad t \rightarrow \infty, \quad \alpha = \frac{t}{n}=\text{finite}
\end{equation}
the quantities ${h_{n,t}(v)}$ simplify:
\begin{align}
\label{eq:4A1}
\begin{split}
    h_{n,t}(v)\to h_\alpha(v) & = \lim_{n \rightarrow \infty, t/n \rightarrow \alpha} \binom{n-1}{t-1}^{-1} \binom{n-v-1}{t-2} \\
    & = \alpha (1-\alpha)^{v-1}.
\end{split}
\end{align}
We shall use the scaled quantities $h_\alpha(v)$ later.

\subsection{Proof of Theorem~\ref{the:ak}}
\label{ssec:ak_proof}

\noindent
We compute~${\bar{r}_k}$ for~${k\geq 1}$ using the recursion relations derived in Sec.~\ref{ssec:recursive_method}. We employ generating functions to transform these recursions into differential equations.

We start with recursion~\eqref{eq:r_outline1} in the case~${t=2}$. Using the relations~${c_{n,2}^k=c_{n,1}^k}$ and~${h_{n,2}=\frac{1}{n-1}}$, the recursion can be rewritten in the form
\begin{align}
\label{cn1:rec}
    & (n-1) c_{n,1}^k = c_{n-1,1}^{k-1} + \sum_{v=2}^{n-1} c_{n-v,1}^{k-1} \cdot c_{v,1}^k.
\end{align}
Making the substitution~${v \rightarrow n-v}$ in the sum, we rewrite \eqref{cn1:rec} as 
\begin{align}
\label{eq:4B1}
    & (n-1) c_{n,1}^k = c_{n-1,1}^{k-1} + \sum_{v=1}^{n-2} c_{v,1}^{k-1} \cdot c_{n-v,1}^k .
\end{align}
Now, we define the following generating function for~${|x|<1}$:
\begin{align}
\label{eq:4B2}
    & f_k(x) = \sum_{v=1}^\infty c_{v,1}^k \cdot x^{v-1}.
\end{align}
To write the recursion~\eqref{eq:4B1} in terms of the generating function, we multiply both sides of~\eqref{eq:4B1} by ${x^{n-2}}$ and sum over all ${n \geq 2}$, which yields the differential equation
\begin{align}
    & \frac{df_k(x)}{dx}  = f_{k-1}(x) + \left( f_k(x) - c_{1,1}^k \right) f_{k-1}(x).
\end{align}
Note that for ${k \geq 1}$, ${c_{1,1}^k = 0}$. This differential equation is separable, and we solve it subject to the initial condition~${f_k(0) = c_{1,1}^k = 0}$ to yield
\begin{align}
\label{eq:4B3}
    & f_k(x) = \exp\!\left[\int_0^x f_{k-1}(x^\prime) dx^\prime \right] - 1.
\end{align}
We get the initial condition in~${k}$ by plugging~${k=0}$ and~${c_{v,1}^0 = 1}$ into Eq.~\eqref{eq:4B2}, which yields~${f_0(x) = \frac{1}{1-x}}$. Iterating the recurrence \eqref{eq:4B3}, we find 
\begin{align}
\label{eq:fk_small_k}
\begin{split}
    f_0(x) &= \frac{1}{1-x}, \\
    f_1(x) &= \frac{x}{1-x}, \\
    f_2(x) &= \frac{e^{-x}}{1-x} - 1, \\
    f_3(x) &= \exp\!\left( -x + \frac{\text{Ei}(1)-\text{Ei}(1-x)}{e} \right) - 1,
\end{split}
\end{align}
where~${\text{Ei}(x) = - \int_{-x}^\infty \frac{e^{-t}}{t} dt}$ is the exponential integral. Next, we analyze \eqref{eq:r_outline1} for general~${t \geq 2}$:
\begin{align}
    \nonumber & c_{n,t}^k = 
    h_{n,t}(1) ~\mathbbm{E}_{t^\prime<t} \left[ c_{n-1,t'}^{k-1} \right] + 
    \sum_{v=2}^{n-t+1} h_{n,t}(v) ~\mathbbm{E}_{t'<t} \left[ c_{n-v,t'}^{k-1} \cdot c_{v,1}^k \right].
\end{align}
In the scaling limit \eqref{scaling}, the quantities ${c_{n,t}^k}$ approach to ${c_k(\alpha)}$. We will also use the function
\begin{align}
    & g_k(\alpha) = \lim_{n \rightarrow \infty, t/n \rightarrow \alpha} \frac{1}{n} \sum_{t^\prime = 1}^t c_{n,t'}^k
    = \int_0^\alpha c_k(\alpha^\prime) d\alpha^\prime.
\end{align}
Note that for~${k=0}$,
\begin{align}
    c_0(\alpha) & = \lim_{n \rightarrow \infty, t/n \rightarrow \alpha} c_{n,t}^0 = \lim_{n \rightarrow \infty, t/n \rightarrow \alpha} 1 = 1,
\end{align}
and thus~${g_0(\alpha) = \alpha}$. Taking the scaling limit \eqref{scaling} and using Eq.~\eqref{eq:4A1}, we recast the recurrence \eqref{eq:r_outline1} into the differential equation
\begin{align}
    & \frac{d}{d\alpha} g_k(\alpha) = g_{k-1} (\alpha) + \sum_{v=2}^\infty (1 - \alpha)^{v-1} g_{k-1}(\alpha) \cdot c_{v,1}^k,
\end{align}
which we simplify using the generating function ${f_k(x)}$ to obtain
\begin{align}
    & \frac{d}{d\alpha} g_k(\alpha) = g_{k-1} (\alpha) + g_{k-1}(\alpha) f_k(1 - \alpha).
\end{align}
Solving this equation yields the following recursion in~${k}$ for ${g_k(\alpha)}$:
\begin{align}
    & g_k(\alpha) = \int_0^\alpha g_{k-1} (\alpha^\prime) \left( 1 + f_k(1 - \alpha^\prime) \right) d \alpha^\prime.
\end{align}
We solve this recurrence iteratively and express ${g_k(\alpha)}$ in terms of the~${f_k}$:
\begin{align}
\label{eq:4B4}
    & g_k(\alpha) = \int\displaylimits_{0<\alpha_1<\cdots<\alpha_k<\alpha} \alpha_1 \prod_{i=1}^k \left( 1 + f_i(1 - \alpha_i) \right) d\alpha_1 \cdots d\alpha_k.
\end{align}

We now express ${b_k(\alpha)}$,  the scaling limit of~${b_{n,t}^k}$, via ${g_k(\alpha)}$. Using \eqref{eq:r_outline1}, we rewrite Eq.~\eqref{eq:r_outline5} as
\begin{align}
\begin{split}
    & b_{n,t}^k = c_{n,t}^k - \frac{1}{n-1} \sum_{t^\prime=1}^{t-1} c_{n-1,t'}^{k-1}
\end{split}
\end{align}
which in the scaling limit becomes $b_k(\alpha) = \frac{d g_k(\alpha)}{d\alpha} - g_{k-1}(\alpha)$ from which
\begin{align}
\label{rk_g}
\begin{split}
    \bar{r}_k & = \displaystyle\lim_{n \rightarrow \infty} \frac{1}{n} \sum_{t=1}^n b_{n,t}^k 
    = \int_0^1 \left( \frac{d}{d\alpha} g_k(\alpha) - g_{k-1}(\alpha) \right)  d \alpha \\
    & = g_k(1) - \int_0^1 g_{k-1}(\alpha) d \alpha.
\end{split}
\end{align}
Using Eq.~\eqref{eq:4B4} for~${g_k(\alpha)}$, we transform \eqref{rk_g} into
\begin{align}
\begin{split}
    \bar{r}_k & = \int\displaylimits_{0<\alpha_1<\cdots<\alpha_k<1}  \alpha_1 \left[ \prod_{i=1}^k \left( 1 + f_i(1 - \alpha_i) \right) - \prod_{i=1}^{k-1} \left( 1 + f_i(1 - \alpha_i) \right)\right] d\alpha_1 \cdots d\alpha_k \\
    & = \int\displaylimits_{0<\alpha_1<\cdots<\alpha_k<1} \alpha_1 f_k(1 - \alpha_k) \prod_{i=1}^{k-1} \left( 1 + f_i(1 - \alpha_i) \right) d\alpha_1 \cdots d\alpha_k.
\end{split}
\end{align}

\subsection{Proof of Theorem \ref{the:pk}}
\label{ssec:pk_proof}

\noindent
Here we compute~${\bar{\ell}_k}$ for~${k \geq 1}$ by solving the recursions derived in Sec.~\ref{ssec:recursive_method}. We rely on the same approach as in the previous section~\ref{ssec:ak_proof}. First, we determine ${c_{n,1}^k}$. Plugging~${t=2}$ into the recursion~\eqref{eq:l_outline1} and using~${c_{n,2}^k=c_{n,1}^k}$ and~${h_{n,2} = \frac{1}{n-1}}$ as before, we obtain a recurrence for ${1 - c_{n,1}^k}$,
\begin{align}
\label{eq:4B7}
    (n-1)(1-c_{n,1}^k) = \sum_{v=1}^{n-1} (1 - c_{n-v,1}^{k-1})(1 - c_{v,1}^k),
\end{align}
which we treat using the generating function
\begin{align}
\label{eq:4B5}
    q_k(x) = \sum_{v\geq 1} (1-c_{v,1}^k) x^{v-1}.
\end{align}
Multiplying \eqref{eq:4B7} by ${x^{n-2}}$ and summing over all ${n}$, we arrive at the differential equation
\begin{align}
    \frac{dq_k(x)}{dx} = q_k(x) q_{k-1}(x) 
\end{align}
which we solve subject to the initial condition $q_k(0) = c_{1,1}^k = 1$ and find
\begin{align}
    q_k(x) = \exp\!\left[\int_0^x q_{k-1}(x^\prime) dx^\prime \right]
\end{align}
for $k\geq 1$. (Note that $q_0(x) = \sum_{v\geq 1} (1-c_{v,1}^0) x^{v-1} = 1$.)

It is not necessary to compute the scaling limit of~${c_{n,t}^k}$, since for large~${n}$, the~${b_{n,t}^k}$ can be expressed in terms of the~${c_{v,1}^k}$, Eq.~\eqref{eq:l_outline3}. In the scaling limit, Eq.~\eqref{eq:l_outline3} becomes
\begin{equation}\begin{aligned}
    b_k(\alpha) & = \lim_{n \rightarrow \infty, t/n \rightarrow \alpha} b_{n,t}^k \\
    & = \lim_{n \rightarrow \infty, t/n \rightarrow \alpha} \sum_{v=1}^{n-t+1} h_{n,t}(v) c_{v,1}^k\\
    & = \sum_{v=1}^\infty \alpha (1 - \alpha)^{v-1} \cdot c_{v,1}^k \\
    & = 1 - \alpha q_k(1 - \alpha)
\end{aligned}\end{equation}
where in the last two steps we used \eqref{eq:4A1} and \eqref{eq:4B5}. 
Finally, ${\bar{\ell}_k}$ is given by 
\begin{equation}
\label{ell-k}
     \bar{\ell}_k  = \displaystyle\lim_{n \rightarrow \infty} \frac{1}{n} \sum_{t=1}^n b_{n,t}^k
     = 1 - \int_0^1 \alpha q_k(1 - \alpha) d \alpha.
\end{equation}

To rewrite this expression in the form given in Theorem~\ref{the:pk}, we introduce a new set of functions
\begin{equation}
s_k(y) = 1 - e^{-y} q_{k-1}(1-e^{-y})
\end{equation}
for~${k \geq 1}$. The recursion relation for the~${s_k(y)}$ is
\begin{align}
\label{eq:4B6}
    s_{k+1}(y) = 1 - \exp\!\left[-\int_0^y s_k(y^\prime) dy^\prime \right].
\end{align}
The base case is~${s_1(y) = 1 - e^{-y}}$, but since defining ${s_0(y) = 1}$ is consistent with \eqref{eq:4B6}, we re-define the latter as the base case. The first four functions~${s_k(y)}$ are
\begin{align}
\label{eq:sk_small_k}
\begin{split}
    & s_0(y) = 1, \\
    & s_1(y) = 1 - e^{-y}, \\
    & s_2(y) = 1 - \exp \left( 1 - y - e^{-y} \right), \\
    & s_3(y) = 1 - \exp\!\left[ -1 - y + \exp \left( 1 - e^{-y} \right) \right]. 
\end{split}
\end{align}
Equation \eqref{ell-k} becomes 
\begin{align}
\label{eq:layer_solution}
    \bar{\ell}_k = 1 - \int_0^1 \left[ 1 - s_k(-\log \alpha) \right] d \alpha = \int_0^\infty e^{-y} s_k(y) dy ,
\end{align}
where in the last step we made the substitution~${\alpha \rightarrow e^{-y}}$. As a consistency check we note that the definition $s_0(y)=1$ is consistent with \eqref{eq:layer_solution} giving the correct answer ${\bar{\ell}_0 = 1}$.

\subsection{Calculations in Corollary~\ref{cor:smallk}}
\label{ssec:rl_small_k}

\noindent
Here we compute ${r_k}$ and the ${\ell_k}$ for $k\leq 2$. By definition, ${\bar{r}_0=\bar{\ell}_0=1}$. Using Theorem~\ref{the:ak} and ${f_i(x)}$ for $i\leq 3$ given in Eq.~\eqref{eq:fk_small_k} we compute the cumulative rank distribution $\bar{r}_k$ for $k=1,2,3$:
\begin{equation}\begin{aligned}
    \bar{r}_1 & = \int\displaylimits_0^1 (1-\alpha_1) d\alpha_1 = \frac{1}{2}\,, \\
    \bar{r}_2 & = \int\displaylimits_{0<\alpha_1<\alpha_2<1} \left( \frac{e^{\alpha_2-1}}{\alpha_2} - 1 \right) d \alpha_1 d \alpha_2
    = \frac{1}{2} - \frac{1}{e}\,, \\
    \bar{r}_3 & = \int\displaylimits_{0 < \alpha_1 < \alpha_2 < \alpha_3 < 1} \left(\exp\! \left[ \alpha_3-1 + \frac{\text{Ei}(1)-\text{Ei}(\alpha_3)}{e} \right] - 1 \right)  \frac{e^{\alpha_2-1}}{\alpha_2} d\alpha_1 d\alpha_2 d\alpha_3 \\
    & = -1 + \frac{2}{e} + \frac{1}{e^2} \int\displaylimits_0^1 \left( e^{2 \alpha_3} - e^{\alpha_3} \right) \exp\! \left[ \frac{\text{Ei}(1) - \text{Ei}(\alpha_3)}{e} \right] d \alpha_3 .
\end{aligned}\end{equation}
Therefore
\begin{equation} \begin{aligned}
     & r_0 = \bar{r}_0 - \bar{r}_1 = \frac{1}{2}, \\
     & r_1 = \bar{r}_1 - \bar{r}_2 = \frac{1}{e}, \\
     & r_2 = \bar{r}_2 - \bar{r}_3 = \frac{3}{2} - \frac{3}{e} - \frac{1}{e^2} \int\displaylimits_0^1 \left( e^{2 \alpha_3} - e^{\alpha_3} \right) \exp\! \left[ \frac{\text{Ei}(1) - \text{Ei}(\alpha_3)}{e} \right] d \alpha_3 \approx 0.124.
\end{aligned} \end{equation}

Substituting Eqs.~\eqref{eq:sk_small_k} into \eqref{eq:layer_solution}, we compute the cumulative layer distribution for $k\leq 3$,
\begin{equation}\begin{aligned}
      \bar{\ell}_0 & = 1, \\
     \bar{\ell}_1 &= \displaystyle\int_0^\infty e^{-y} \left( 1 - e^{-y} \right) dy = \frac{1}{2}\,, \\
     \bar{\ell}_2 & = \displaystyle\int_0^\infty e^{-y} \left[1 - \exp \left( 1 - y - e^{-y} \right) \right] dy = 3 - e, \\
     \bar{\ell}_3 & = \displaystyle\int_0^\infty e^{-y} \left( 1 - \exp\!\left[ -1 - y + \exp \left( 1 - e^{-y} \right) \right] \right) dy = 1 - \frac{1}{e} \displaystyle\int_0^1 u \exp \left( e^{1 - u} \right) du,
\end{aligned}\end{equation}
from which we recover the layer distribution for $k\leq 2$:
\begin{equation} \begin{aligned}
    & \ell_0 = \bar{\ell}_0 - \bar{\ell}_1 = \frac{1}{2}, \\
    & \ell_1 = \bar{\ell}_1 - \bar{\ell}_2 = e - \frac{5}{2}, \\
    & \ell_2 = \bar{\ell}_2 - \bar{\ell}_3 = 2 - e + \frac{1}{e} \displaystyle\int_0^1 u \exp \left( e^{1 - u} \right) du \approx 0.114.
\end{aligned} \end{equation}

\subsection{Proof of Theorem~\ref{the:rank_dist_tail}}
\label{ssec:rank_tail_proof}

\noindent
In this section, we prove Eq.~\eqref{eq:rank_tail}. Recall the result of Theorem~\ref{the:ak}:
\begin{equation}\begin{aligned}
    &\bar{r}_k  = \int\displaylimits_{0 < \alpha_1 < \dots < \alpha_k < 1} \alpha_1 f_k(1-\alpha_k) \prod\displaylimits_{i=1}^{k-1} \left( 1 + f_i(1-\alpha_i) \right) d \alpha_1 \dots d \alpha_k \\
   & f_0(x)  = \frac{1}{1-x}, \quad
    f_i(x)  = \exp \left(\int_0^x f_{i-1}(x^\prime) dx^\prime \right) - 1 \quad \forall i \geq 1 ,
\end{aligned}\end{equation}
where the~${f_i(x)}$ are defined when $0\leq x <1$. First, we show by induction that for all~${k \geq 4}$,
\begin{align}
\label{Efk_largeK}
     \frac{x^k}{k!} \leq f_k(x) \leq \frac{x^{k-2}}{(k-2)!} .
\end{align}

The lower bound is valid for all $k\geq 0$. Indeed, $f_0(x) = \frac{1}{1-x} \geq 1$, and assuming the bound to hold in the rank $k$, we get 
\begin{align}
     f_{k+1}(x) = \exp \left( \int_0^x f_k(x^\prime) dx^\prime \right) - 1
    \geq \exp \left( \frac{x^{k+1}}{(k+1)!} \right) - 1
    \geq \frac{x^{k+1}}{(k+1)!} 
\end{align}
completing the proof by induction.

We now prove by induction a stronger assertion, ${f_k(x) \leq \frac{x^{k-2}}{(k-2)!} - \frac{x^{k-2}}{(k-1)!}}$ for $k\geq 4$, which implies the right inequality in \eqref{Efk_largeK}. We first check 
\begin{align}
\label{f4}
    f_4(x) \leq \frac{x^2}{2!} - \frac{x^2}{3!}
\end{align}
for ${x \in [0,1]}$. This is easy to see by plotting $f_4(x)$. The small $x$ expansion $f_4(x)=\frac{x^4}{24}+\frac{x^5}{60}+\frac{x^6}{80}+\ldots$ also supports the inequality \eqref{f4}. One could wonder about the applicability of \eqref{f4} at $x=1$ where $f_4$ is singular. Indeed, the explicit formulas \eqref{eq:fk_small_k} show that the functions $f_k$ with $k\leq 2$ diverge as $(1-x)^{-1}$, while the divergence of $f_3$ is weaker: $f_3\sim (1-x)^{-1/e}$. The function $f_4$ is singular at $x=1$, e.g., its derivatives diverge, but the limiting value is finite, $f_4(1)\approx 0.210719$, so \eqref{f4} remains valid. 

To prove the inductive step, three additional inequalities are necessary. First, using elementary calculus one can show that for~${y \in [0,1]}$,
\begin{align}
    \label{Eineq_H1}
    & e^y - 1 \leq y + y^2.
\end{align}
Additionally, for~${k \geq 4}$ we have
\begin{align}
    \label{Eineq_H2}
    & \frac{x^{k+1}}{(k+1)!} - \frac{x^{k+1}}{(k+1)(k+1)!}
    \leq \frac{1}{(k+1)!} - \frac{1}{(k+1)(k+1)!}
    \leq \frac{1}{150}
    < 1.
\end{align}
Finally, one can also show that for~${k \geq 4}$,
\begin{align}
    \label{Eineq_H3}
    & \left( \frac{1}{k} - \frac{1}{k-1} \right) + \frac{1}{(k-1)!} \left( 1 - \frac{1}{k-1} \right)^2 < 0.
\end{align}
Using Eqs.~\eqref{Eineq_H1}--\eqref{Eineq_H3}, we prove the inductive step. Assuming~${f_k(x) \leq \frac{x^{k-2}}{(k-2)!} - \frac{x^{k-2}}{(k-1)!}}$, we have
\begin{equation}\begin{aligned}
    & f_{k+1}(x) = \exp \left( \int_0^x f_k(x^\prime) dx^\prime \right) - 1 \\
    & \leq \exp \left( \frac{x^{k-1}}{(k-1)!} - \frac{x^{k-1}}{(k-1)(k-1)!} \right) - 1 \\
    & \leq \left( \frac{x^{k-1}}{(k-1)!} - \frac{x^{k-1}}{(k-1)(k-1)!} \right) + \left( \frac{x^{k-1}}{(k-1)!} - \frac{x^{k-1}}{(k-1)(k-1)!} \right)^2 \\
    & \leq \left( \frac{x^{k-1}}{(k-1)!} - \frac{x^{k-1}}{k!} \right) + \frac{x^{k-1}}{(k-1)!} \left( \frac{1}{k} - \frac{1}{k-1} \right) + \left( \frac{x^{k-1}}{(k-1)!} - \frac{x^{k-1}}{(k-1)(k-1)!} \right)^2 \\
    & \leq \left( \frac{x^{k-1}}{(k-1)!} - \frac{x^{k-1}}{k!} \right) + \frac{x^{k-1}}{(k-1)!} \left[ \left( \frac{1}{k} - \frac{1}{k-1} \right) + \frac{1}{(k-1)!} \left( 1 - \frac{1}{k-1} \right)^2 \right] \\
    & \leq \frac{x^{k-1}}{(k-1)!} - \frac{x^{k-1}}{k!}.
\end{aligned}\end{equation}
This completes the induction.

Next, we derive the asymptotic value of~${r_k}$. For all~${x \in [0,1)}$, we have by Eq.~\eqref{Efk_largeK} that
\begin{align}
\label{eq:rank_tail1}
    1 \leq \prod_{i=4}^{k-1} \left( 1 + f_i(x_i) \right)
    \leq \prod_{i=4}^{k-1} \left( 1 + \frac{1}{(i-2)!} \right)
    < \prod_{i=4}^{k-1} \exp \left( \frac{1}{(i-2)!} \right)
    < e^{e-2}.
\end{align}

From the expression for~${\bar{r}_k}$, we derive the following expression for ${r_k = \bar{r}_k - \bar{r}_{k+1}}$:
\begin{equation}\begin{aligned}
\label{eq:rank_tail2}
    r_k  =& \int\displaylimits_{0 < \alpha_1 < \dots < \alpha_{k+1} < 1} \alpha_1 \left[ \frac{f_k(1-\alpha_k)}{1-\alpha_k} - f_{k+1}(1-\alpha_{k+1}) \left( 1 + f_k(1-\alpha_k) \right)\right] \nonumber\\ &\times \prod\displaylimits_{i=1}^{k-1} \left( 1 + f_i(1-\alpha_i) \right) d \alpha_1 \dots d \alpha_{k+1},
\end{aligned}\end{equation}

To establish the upper bound we rely on \eqref{eq:rank_tail1} and write
\begin{equation}\begin{aligned}
    r_k & \leq e^{e-2} \int\displaylimits_{0 < \alpha_1 < \dots < \alpha_{k+1} < 1} \alpha_1 \frac{f_k(1-\alpha_k)}{1-\alpha_k} \prod\displaylimits_{i=1}^3 \left( 1 + f_i(1-\alpha_i) \right) d \alpha_1 \dots d \alpha_{k+1}.
\end{aligned}\end{equation}
which we further amend using the upper bound in \eqref{Efk_largeK} to yield 
\begin{equation}\begin{aligned}
    r_k & \leq e^{e-2} \int\displaylimits_{0 < \alpha_1 < \dots < \alpha_{k+1} < 1} \alpha_1 \frac{(1-\alpha_k)^{k-3}}{(k-2)!} \prod\displaylimits_{i=1}^3 \left( 1 + f_i(1-\alpha_i) \right) d \alpha_1 \dots d \alpha_{k+1} \\
    & \leq \frac{e^{e-2}}{(2k-5)!} \int\displaylimits_{0 < \alpha_1 < \dots < \alpha_3 < 1} \alpha_1 \prod\displaylimits_{i=1}^3 \left( 1 + f_i(1-\alpha_i) \right) d \alpha_1 \dots d \alpha_3 \\
    & = \frac{e^{e-2}}{(2k-5)!} \left( \bar{r}_3 + 1 - \frac{2}{e} \right)\\
    & \leq \exp\!\left[ -2k \log k + (2 - 2\log2) k + O(\log k) \right]
\end{aligned}\end{equation}
where in the last step we used Stirling's formula.

To prove the lower bound we use Eq.~\eqref{eq:rank_tail2} for~${r_k}$ and the fact that~${f_i(x)}$ is positive:
\begin{equation}\begin{aligned}
\label{eq:rank_tail3}
    r_k & \geq \int\displaylimits_{0 < \alpha_1 < \dots < \alpha_{k+1} < 1} \alpha_1 \left[ \frac{f_k(1-\alpha_k)}{1-\alpha_k} - f_{k+1}(1-\alpha_{k+1}) \left( 1 + f_k(1-\alpha_k) \right)\right] d \alpha_1 \dots d \alpha_{k+1} \\
    & = \int\displaylimits_{0 < \alpha_1 < \dots < \alpha_k < 1} \alpha_1 \left[ f_k(1-\alpha_k) - \left( 1 + f_k(1-\alpha_k) \right) \int_{\alpha_k}^1 f_{k+1}(1-\alpha_{k+1}) d\alpha_{k+1} \right] d \alpha_1 \dots d \alpha_k .
\end{aligned}\end{equation}
It is straightforward to prove by induction that the~${f_i(x)}$ are all increasing. Hence
\begin{equation}\begin{aligned}
    \int_{\alpha_k}^1 f_{k+1}(1-\alpha_{k+1}) d\alpha_{k+1} & \leq (1-\alpha_k) f_{k+1}(1-\alpha_k) \\
    & = (1-\alpha_k) \left( \exp \left( \int_0^{1-\alpha_k} f_k(x) dx \right) - 1 \right) \\
    & \leq (1-\alpha_k) \left( \exp \left( (1-\alpha_k) f_k(1-\alpha_k) \right) - 1 \right) \\
    & = (1-\alpha_k)^2 f_k(1-\alpha_k) + O(f_k(1-\alpha_k)^2).
\end{aligned}\end{equation}
Plugging this into~\eqref{eq:rank_tail3}, we obtain
\begin{equation}\begin{aligned}
    r_k & \geq \int\displaylimits_{0 < \alpha_1 < \dots < \alpha_k < 1} \alpha_1 \left[ (1 - (1-\alpha_k)^2) f_k(1-\alpha_k) + O(f_k(1-\alpha_k)^2) \right] d \alpha_1 \dots d \alpha_k.
\end{aligned}\end{equation}

Using~\eqref{Efk_largeK}, we can further lower bound this expression by
\begin{equation}\begin{aligned}
    r_k & \geq \int\displaylimits_{0 < \alpha_1 < \dots < \alpha_k < 1} \alpha_1 \left[ (1 - (1-\alpha_k)^2) \frac{(1-\alpha_k)^k}{k!} + O\left(\frac{1}{((k-2)!)^2}\right) \right] d \alpha_1 \dots d \alpha_k .
\end{aligned}\end{equation}
Evaluating the integral yields
\begin{equation}\begin{aligned}
    r_k & \geq \frac{(k+1)(3k+4)}{(2k+3)!} + O\left(\frac{1}{(k+1)!((k-2)!)^2}\right) \\
    & \geq  \exp\!\left[ -2k \log k + (2-2\log2)k + O(\log k) \right]
 \end{aligned}\end{equation}
where in the last step we used  Stirling's formula. This completes the proof of the asymptotic Eq.~\eqref{eq:rank_tail}.

\subsection{Proof of Theorem~\ref{the:layer_dist_tail}}
\label{ssec:layer_tail_proof}

\noindent
Here we prove exponential bounds for  the cumulative layer distribution
\begin{equation}
\label{ell-k:bounds}
    e^{-k} \lesssim \bar{\ell}_k \lesssim e^{-k/e}.
\end{equation}
Theorem~\ref{the:pk} states that $\bar{\ell}_k = \int_0^\infty dy~ e^{-y} s_k(y)$ with $s_0(y)  = 1$ and $s_i(y)$ for $i>0$ determined recursively via 
\begin{equation*}
s_i(y)  = 1 - \exp\!\left[ -\int_0^y dy^\prime~ s_{i-1}(y^\prime) \right]
\end{equation*}
when $i\geq 1$. All functions ${s_i(y)}$ are defined on the interval~${y \in [0,\infty)}$ and satisfy ${0 < s_i(y) \leq 1}$. 

To prove the lower bound in Eq.~\eqref{eq:layer_bounds}, we take $\epsilon \in (0,1)$ and note that~${s_1(y)}$ has the property
\begin{align}
    & s_1(y) = 1 - e^{-y} > 1 - \epsilon \quad \text{for}~ y > - \log \epsilon .
\end{align}
We will prove by induction that each~${s_k(y)}$ has a similar property, defined below. Assume that for some~${k \geq 1}$ and~${Y > 0}$,
\begin{align}
    & s_k(y) > 1 - \epsilon \quad \text{for}~ y > Y .
\end{align}
Then for~${y > Y + \frac{- \log \epsilon}{1 - \epsilon}}$, 
\begin{equation}\begin{aligned}
    s_{k+1}(y) & = 1 - e^{- \int_0^y s_k(y^\prime) dy^\prime} \\
    & \geq 1 - e^{- \int_Y^{Y+\frac{- \log \epsilon}{1 - \epsilon}} (1-\epsilon) dy^\prime} \\
    & \geq 1 - e^{- (1 - \epsilon) \frac{- \log \epsilon}{1 - \epsilon}} \\
    & = 1 - \epsilon .
\end{aligned}\end{equation}
By induction we have, for any~${k \geq 1}$,
\begin{align}
    & s_k(y) > 1 - \epsilon \quad \text{for}~ y > - \log \epsilon \left( 1 + \frac{k}{1 - \epsilon} \right) .
\end{align}
Therefore, for any~${k}$ we have 
\begin{align}
    & \bar{\ell}_k = \int_0^\infty s_k(y) e^{-y} dy \geq \int_{- \log \epsilon \left( 1 + \frac{k}{1 - \epsilon} \right)}^\infty (1 - \epsilon) e^{-y} dy = (1 - \epsilon) \epsilon^{1 + \frac{k}{1 - \epsilon}} ,
\end{align}
from which we find an asymptotic lower bound on~${\bar{\ell}_k}$: ${\bar{\ell}_k \gtrsim \epsilon^{\frac{k}{1 - \epsilon}}}$. It remains to choose the optimal value of~${\epsilon}$ which will give the tightest lower bound. Setting~${\epsilon = 1/2}$, we obtain~${\bar{\ell}_k \gtrsim 4^{-k}}$. However, setting~${\epsilon}$ arbitrarily close to~${1}$, we get a better lower bound: ${\bar{\ell}_k \gtrsim e^{-k}}$.

We can use a similar technique to place an asymptotic upper bound on~${\bar{\ell}_k}$. We notice that the inequality $s_1(y) = 1 - e^{-y} < y$ is valid for  $y>0$, and assume that it is a special case of the inequalities
\begin{align}
\label{sj<}
    s_k(y) < \frac{y^{k+1}}{(k+1)!}
\end{align}
valid for all $k\geq 1$ and all ${y > 0}$. A straightforward computation 
\begin{equation}\begin{aligned}
    s_{k+1}(y) & = 1 - e^{- \int_0^y s_k(y^\prime) dy^\prime} \\
    & < 1 - e^{- \int_0^y \frac{(y^\prime)^{k+1}}{(k+1)!} dy^\prime} \\
    & = 1 - e^{- \frac{y^{k+2}}{(k+2)!}}  < \frac{y^{k+2}}{(k+2)!}
\end{aligned}\end{equation}
proves \eqref{sj<} by induction.

We can strengthen the result~\eqref{sj<} using the trivial upper bound ${s_k(y) < 1}$ for all~${k}$, which is a tighter bound when~${y>[(k+1)!]^{1/(k+1)}}$. Putting both bounds together yields
\begin{equation}\begin{aligned}
    s_k(y) < \min \left[ \frac{y^{k+1}}{(k+1)!}, 1 \right] .
\end{aligned}\end{equation}
Setting ${\gamma(k) = [(k+1)!]^{1/(k+1)}}$, we have 
\begin{equation}\begin{aligned}
\label{ell-k:upper}
    \bar{\ell}_k & = \int_0^\infty s_k(y) e^{-y} dy \\
    & < \int_0^{\gamma(k)} \frac{y^{k+1}}{(k+1)!} e^{-y} dy + \int_{\gamma(k)}^\infty 1 \cdot e^{-y} dy \\
    & = \left[ -\frac{\gamma(k)^{k+1}}{(k+1)!} e^{-\gamma(k)} -\frac{\gamma(k)^k}{k!} e^{-\gamma(k)} - \dots - e^{-\gamma(k)} + 1 \right] + e^{-\gamma(k)} \\
    & = e^{-\gamma(k)} \left[ 1 + \frac{(k+1)! \gamma(k)}{(k+2)!} + \frac{(k+1)! \gamma(k)^2}{(k+3)!} + \dots \right] \\
    & < e^{-\gamma(k)} \left[ 1 + \frac{\gamma(k)}{k+2} + \frac{\gamma(k)^2}{(k+2)^2} + \dots \right] 
    = \frac{e^{-\gamma(k)}}{1 - \frac{\gamma(k)}{k+2}} \,.
\end{aligned}\end{equation}
Stirling's formula gives ${\gamma(k) = [(k+1)!]^{1/(k+1)}}\simeq k/e$ which we insert into \eqref{ell-k:upper} and arrive at the upper 
bound in \eqref{ell-k:bounds}.

\subsection{Proof of Theorem \ref{the:Rnk_variance}}
\label{ssec:variance_proof}

\noindent
Theorem~\ref{the:Rnk_variance} asserts that the number of nodes of rank ${k}$ is strongly self-averaging:
\begin{equation}\begin{aligned}
    & \text{Var} \left[ \frac{R_{n,k}}{n} \right] \lesssim \frac{1}{n}\, .
\end{aligned}\end{equation}
Here ${R_{n,k}}$ is a random variable which denotes the number of nodes of rank~${k}$ in an RRT of size~${n}$. The theorem also holds when~${R_{n,k}}$ is replaced by~${L_{n,k}}$, the number of nodes in layer $k$; all steps in the proof below can be carried out for $L_{n,k}$.

The main tool in the proof will be RRT \emph{splitting,} where the edge between nodes~${1}$ and~${2}$ is deleted to obtain 
two smaller RRTs of sizes~${i}$ and~${n-i}$, where~${i}$ is uniformly distributed from~${1}$ to~${n-1}$. However, RRT splitting 
can potentially change the ranks of nodes that are near the split. Let~${X_{n,k}}$ denote the error due to splitting: 
\begin{align}
    & R_{n,k} = R_{i,k} + R_{n-i,k} + X_{n,k}.
\end{align}

A node of rank/layer number~${k}$ before splitting can have a different rank/layer number after splitting  only if its distance to node~${1}$ (the root) is~${\leq k}$. Let~${M_{n,r}}$ denote the number of nodes a distance~${r}$ from node~${1}$. 
The number of RRTs with~${M_{n,r-1} \leq m}$ and~${M_{n,r} = i}$ has the upper bound
\begin{align}
    & m^i \left( \sum_{1 \leq \sigma_1 < \sigma_2 < \dots < \sigma_{n-i-1} \leq n-1 } \prod_j \sigma_j \right) 
\end{align}
since there are at most~${m}$ places to attach each of the~${i}$ nodes a distance~${r}$ away from the root, and they can be 
placed at any times from~${1}$ to~${n}$ during the RRT construction.

For any~${c > 0}$, we will show via induction 
on~${r}$ that 
\begin{align}
\label{PMnr}
    P (M_{n,r} > c (\log n)^{r+1}) \lesssim n^{-c/e^{2r}} .
\end{align}
Indeed, ${M_{n,0} = 1}$, and hence  $P (M_{n,0} > c \log n) = 0 < n^{-c}$, i.e., \eqref{PMnr} is obeyed when $r=0$. 
Assuming the hypothesis \eqref{PMnr} for~${r-1}$, we have 
\begin{equation}\begin{aligned}
& P (M_{n,r} > c (\log n)^{r+1})  \leq P[M_{n,r-1} > (c/e^2) (\log n)^r] + \Pi \\
& \Pi = P[M_{n,r} > c (\log n)^{r+1} ~\text{and}~ M_{n,r-1} \leq (c/e^2) (\log n)^r]
\end{aligned}\end{equation}
Using \eqref{PMnr}, we get the bound for the first term 
\begin{equation}
\label{ineq:1}
P[M_{n,r-1} > (c/e^2) (\log n)^r] \lesssim n^{-(c/e^2)/e^{2(r-1)}} = n^{-c/e^{2r}} 
\end{equation}
and then the bound for the second term 
\begin{equation}\begin{aligned}
\label{ineq:2}
\Pi    & \leq \frac{1}{(n-1)!} \sum_{i > c (\log n)^{r+1}} \left( (c/e^2) (\log n)^r \right)^i 
    \left( \sum_{1 \leq \sigma_1 < \sigma_2 < \dots < \sigma_{n-i-1} \leq n-1 } \prod_j \sigma_j \right) \\
    & \leq  \frac{1}{(n-1)!} \sum_{i > c (\log n)^{r+1}} \left( (c/e^2) (\log n)^r \right)^i 
    \left( \frac{(n-1)!}{i!} \left( 1 + \frac{1}{2} + \dots + \frac{1}{n-1} \right)^i \right) \\
    & \leq  \frac{1}{(n-1)!} \sum_{i > c (\log n)^{r+1}} 
    \left( (c/e^2) (\log n)^r \right)^i \left( \frac{(n-1)!}{i!} (\log n)^i \right) \\
    & =  \sum_{i > c (\log n)^{r+1}} \frac{\left( (c/e^2) (\log n)^{r+1} \right)^i}{i!} < \sum_{i > c (\log n)^{r+1}} \frac{\left( (c/e^2) (\log n)^{r+1} \right)^i}{(i/e)^i} \\
    & < \sum_{i > c (\log n)^{r+1}} e^{-i}  = \frac{e^{-c (\log n)^{r+1}}}{1 - 1/e}\,. 
\end{aligned}\end{equation}
In the penultimate step, we used $i!>(i/e)^i$ following from the Stirling formula. Comparing \eqref{ineq:1} and \eqref{ineq:2}, we conclude that $\Pi$ decays faster than any power law for $r\geq 1$, so it is negligible, and we arrive at \eqref{PMnr} thereby completing the induction.

Setting~${c = 2 e^{2r}}$, we obtain~${P (M_{n,r} > 2e^{2r} (\log n)^{r+1}) \lesssim n^{-2}}$. 
This result leads to the following upper bound on the variance of~${X_{n,k}}$: 
\begin{equation}\begin{aligned}
\label{eq:Xvar}
    \text{Var} [X_{n,k}] & \leq \langle X_{n,k}^2 \rangle \\
    & \leq \mathbb{E} \left[ M_{n,k}^2 \right] \\
    & \leq 4e^{4k} (\log n)^{2k+2} P \left( M_{n,k} \leq 2e^{2k} (\log n)^{k+1} \right) + n^2 P \left( M_{n,k} > 2e^{2k} (\log n)^{k+1} \right) \\
    & \lesssim 4e^{4k} (\log n)^{2k+2} + n^2 \cdot n^{-2} \sim (\log n)^{2k+2} .
\end{aligned}\end{equation}

We now split an RRT in two by deleting the edge between nodes~${1}$ and~${2}$. The resulting RRTs have sizes ${i}$ and ${n-i}$, where~${i}$ is uniformly distributed on~${1,2,\dots,n-1}$. Applying the inequality 
\begin{align}
    & \text{Var} [A+B] \leq \text{Var} [A] + \text{Var} [B] + 2 \sqrt{\text{Var} [A] \text{Var} [B]} ,
\end{align}
we have 
\begin{equation}\begin{aligned}
\label{long}
    \text{Var} [R_{n,k}] 
    & = \frac{1}{n-1} \sum_{i=1}^{n-1} \text{Var} \left[ R_{i,k} + R_{n-i,k} + X_{n,k} \right] \\
    & \leq \frac{1}{n-1} \sum_{i=1}^{n-1} \left[V_{i,n;k} + \text{Var} \left[ X_{n,k} \right]
    + 2 \sqrt{V_{i,n;k} \text{Var} \left[ X_{n,k} \right]} \right] \\
    & = \text{Var} \left[ X_{n,k} \right] + \frac{1}{n-1} \sum_{i=1}^{n-1} V_{i,n;k} 
    + \frac{2}{n-1} \sqrt{\text{Var} \left[ X_{n,k} \right] }\sum_{i=1}^{n-1} \sqrt{V_{i,n;k}}
\end{aligned}\end{equation}
where we used the shorthand notation $V_{i,n;k}=\text{Var} \left[ R_{i,k} + R_{n-i,k} \right]$. Using~\eqref{eq:Xvar}, we now define~${A>0}$ such that~${\text{Var} [X_{n,k}] \leq A (\log n)^{2k+2}}$. We obtain
\begin{equation}\begin{aligned}
    \text{Var} [R_{n,k}]
    & \leq A (\log n)^{2k+2} + \frac{1}{n-1} \sum_{i=1}^{n-1} V_{i,n;k}
    + \frac{2 \sqrt{A} (\log n)^{k+1}}{n-1} \sum_{i=1}^{n-1} \sqrt{V_{i,n;k}}
\end{aligned}\end{equation}
Since~${R_{i,k}}$ and~${R_{n-i,k}}$ are independent, we have~$V_{i,n;k}={\text{Var} \left[ R_{i,k} + R_{n-i,k} \right] = \text{Var} \left[ R_{i,k} \right] + \text{Var} \left[ R_{n-i,k} \right]}$. 
Plugging this into \eqref{long}  and then applying Jensen's inequality to the last term, we have
\begin{equation}\begin{aligned}
    \text{Var} [R_{n,k}] 
    & \leq A (\log n)^{2k+2} 
    + \frac{2}{n-1} \sum_{i=1}^{n-1} \text{Var} \left[ R_{i,k} \right]
    + \frac{2 \sqrt{A} (\log n)^{k+1}}{n-1} \sum_{i=1}^{n-1} \sqrt{\text{Var} \left[ R_{i,k} \right] + \text{Var} \left[ R_{n-i,k} \right]} \\
    & \leq A (\log n)^{2k+2} + \frac{2}{n-1} \sum_{i=1}^{n-1} \text{Var} \left[ R_{i,k} \right] 
    + \frac{2 \sqrt{A} (\log n)^{k+1}}{\sqrt{n-1}} \sqrt{2 \sum_{i=1}^{n-1} \text{Var} \left[ R_{i,k} \right]}.
\end{aligned}\end{equation}

Define~${a_n = \frac{1}{n^2} \sum_{i=1}^n \text{Var} \left[ R_{i,k} \right]}$ for all~${n \geq 1}$. Rewriting the 
inequality above in terms of~${a_n}$ gives
\begin{equation}\begin{aligned}
    & n^2 a_n - (n-1)^2 a_{n-1} \leq A (\log n)^{2k+2} + 2 (n-1) a_{n-1} + 2 \sqrt{2 A} (\log n)^{k+1} \sqrt{n-1} \sqrt{a_{n-1}} \\
    & \Rightarrow a_n \leq \left( 1 - \frac{1}{n^2} \right) a_{n-1} + 2 \sqrt{2 A} \frac{(\log n)^{k+1} \sqrt{n-1}}{n^2} \sqrt{a_{n-1}} + A \frac{(\log n)^{2k+2}}{n^2} \\
    & \Rightarrow a_n \leq a_{n-1} + 2 \sqrt{2 A} \frac{(\log n)^{k+1} \sqrt{n-1}}{n^2} \sqrt{a_{n-1}} + A \frac{(\log n)^{2k+2}}{n^2}.
\end{aligned}\end{equation}

The sequence~${a_n}$ is upper bounded by the sequence~${a^\prime_n}$, where~${a^\prime_1 = a_1}$ and 
\begin{align}
    & a^\prime_n = a^\prime_{n-1} + 2 \sqrt{2 A} \frac{(\log n)^{k+1} \sqrt{n-1}}{n^2} \sqrt{a^\prime_{n-1}} + A \frac{(\log n)^{2k+2}}{n^2}.
\end{align}

Note that~${a^\prime_n}$ is increasing. If it diverges, then there exists~${M}$ such that~${a^\prime_n > 1}$ for~${\forall n \geq M}$. 
Then for~${n \geq M+1}$ we have 
\begin{equation}\begin{aligned}
    & a^\prime_n \leq a^\prime_{n-1} + 2 \sqrt{2 A} \frac{(\log n)^{k+1} \sqrt{n-1}}{n^2} \sqrt{a^\prime_{n-1}} + A \frac{(\log n)^{2k+2}}{n^2} \sqrt{a^\prime_{n-1}} \\
    & \Rightarrow a^\prime_n \leq \left( \sqrt{a^\prime_{n-1}} + \sqrt{2 A} \frac{(\log n)^{k+1} \sqrt{n-1}}{n^2} + \frac{A}{2} \frac{(\log n)^{2k+2}}{n^2} \right)^2 - \left( \frac{A}{2} \frac{(\log n)^{2k+2}}{n^2} \right)^2 \\
    & \Rightarrow \sqrt{a^\prime_n} \leq \sqrt{a^\prime_{n-1}} + \sqrt{2 A} \frac{(\log n)^{k+1} \sqrt{n-1}}{n^2} + \frac{A}{2} \frac{(\log n)^{2k+2}}{n^2}.
\end{aligned}\end{equation}

Iterating gives
\begin{align}
    & \sqrt{a^\prime_n} \leq \sqrt{a^\prime_M} + \sum_{n=M+1}^\infty \left( \sqrt{2 A} \frac{(\log n)^{k+1} \sqrt{n-1}}{n^2} + \frac{A}{2} \frac{(\log n)^{2k+2}}{n^2} \right).
\end{align}

Since the sum on the right converges, the sequence ${a^\prime_n}$ converges, so~${a_n}$ has an upper bound. Then 
\begin{align}
    & \sum_{i=1}^n \text{Var} \left[ R_{i,k} \right] \lesssim n^2 
    \Rightarrow \text{Var} \left[ \frac{R_{n,k}}{n} \right] \lesssim \frac{1}{n}.
\end{align}

\subsection{Proof of Theorem \ref{the:Rnk_expectation}}
\label{ssec:Rnk_expectation}

\noindent
Theorem \ref{the:Rnk_expectation} asserts that the expected value of the fraction of nodes with rank~${k}$ converges and quickly approaches its limiting value:
\begin{equation}
     \Bigg\lvert \Bigg\langle \frac{R_{n,k}}{n} \Bigg\rangle - r_k \Bigg\rvert \lesssim \frac{(\log n)^{k+1}}{n} .
\end{equation}
This theorem also holds for~${L_{n,k}}$ and the corresponding limiting values~${\ell_k}$, and all steps below are equally valid for~${L_{n,k}}$. 

From Eq.~\eqref{eq:Xvar}, we have 
\begin{align}
\label{Rnk}
    \left[ \langle R_{n,k} \rangle - \frac{2}{n-1} \sum_{i=1}^{n-1} \langle R_{i,k} \rangle \right]^2 
    = \langle X_{n,k} \rangle^2 
    \leq \langle X_{n,k}^2 \rangle 
    \leq A (\log n)^{2k+2}
\end{align}
for~${A>0}$ sufficiently large. Now, we define the sequence~${f_n = \langle R_{n,k} \rangle / n}$ and the sequence ${g_n = f_n-f_{n-1}}$ for ${n \geq 1}$, with ${g_0 = f_0}$. 
Rewriting Eq.~\eqref{Rnk} in terms of~${g_n}$ gives 
\begin{equation}
\left[ n \sum_{i=1}^n g_i - \frac{2}{n-1} \sum_{i=1}^{n-1} i \sum_{j=1}^i g_j \right]^2 \leq A (\log n)^{2k+2}
\end{equation}
from which we deduce 
\begin{align}
\label{eq:T5intermediate}
 \Bigg\lvert \sum_{i=1}^n \frac{i(i-1)}{(n+1)n(n-1)} g_i \Bigg\rvert \leq \sqrt{A} \frac{(\log n)^{k+1}}{(n+1)n}.
\end{align}

Summing this inequality over all~${n \geq M}$ and applying the triangle inequality to the LHS gives 
\begin{align}
\begin{split}
    & \Bigg\lvert \sum_{n=M}^\infty \sum_{i=2}^n \frac{i(i-1)}{(n+1)n(n-1)} g_i \Bigg\rvert \leq \sqrt{A} \sum_{n=M}^\infty \frac{(\log n)^{k+1}}{(n+1)n}, \\
    & \Rightarrow \Bigg\lvert \frac{1}{2} \sum_{i=2}^{M-1} \frac{i(i-1)}{M(M-1)} g_i + \frac{1}{2} \sum_{i=M}^\infty g_i \Bigg\rvert \leq \sqrt{A} \sum_{n=M}^\infty \frac{(\log n)^{k+1}}{(n+1)n}, \\
     & \Rightarrow \Bigg\lvert \sum_{i=M}^\infty g_i \Bigg\rvert \leq 2 \sqrt{A} \sum_{n=M}^\infty \frac{(\log n)^{k+1}}{(n+1)n} + \Bigg\lvert \sum_{i=2}^{M-1} \frac{i(i-1)}{M(M-1)} g_i \Bigg\rvert.
\end{split}
\end{align}

Applying the last line in Eqs.~\eqref{eq:T5intermediate} to the second term on the RHS, we have
\begin{align}
\begin{split}
    & \Bigg\lvert \sum_{i=M}^\infty g_i \Bigg\rvert \leq 2 \sqrt{A} \sum_{n=M}^\infty \frac{(\log n)^{k+1}}{(n+1)n} + \sqrt{A} \frac{(M-2)(\log n)^{k+1}}{M(M-1)}, \\
    & \Rightarrow \Bigg\lvert \sum_{i=M}^\infty g_i \Bigg\rvert \lesssim \frac{(\log M)^{k+1}}{M} ,
\end{split}
\end{align}
and replacing~${M}$ with~${n+1}$ and~${\sum_{i=n+1}^\infty g_i}$ with~${r_k-f_n}$, we get 
\begin{align}
    & |f_n - r_k| \lesssim \frac{(\log (n+1))^{k+1}}{n+1} \sim \frac{(\log n)^{k+1}}{n} .
\end{align}

\subsection{Proof of Theorem \ref{the:finite_n}}
\label{ssec:finite_n}

\noindent
From Eq.~\eqref{eq:r_outline1}, we find 
\begin{align}
\label{cn1}
    c_{n,1}^1 & = \frac{1}{n-1} + \frac{1}{n-1} \sum_{v=2}^{n-1} c_{v,1}^1
\end{align}
for~${n \geq 2}$. The base case is ${c_{1,1}^1 = 0}$. The next term ${c_{2,1}^1 = 1}$ follows from \eqref{cn1}. It is then easy to verify that ${c_{n,1}^1 = 1}$ for all $n \geq 2$. Specializing Eq.~\eqref{eq:r_outline1} to ${t \geq 2}$ we obtain 
\begin{equation}\begin{aligned}
    c_{n,t}^1 & = \frac{1}{n-1} \sum_{t^\prime=1}^{t-1} c_{n-1,t'}^0 + \frac{1}{t-1} \binom{n-1}{t-1} ^{-1} \times
    \sum_{v=2}^{n-t+1} \binom{n-v-1}{t-2} \sum_{t^\prime=1}^{t-1} c_{n-v,t'}^0 \times c_{v,1}^1 \\
    & = \frac{t-1}{n-1} + 
    \binom{n-1}{t-1} ^{-1} \times
    \sum_{v=2}^{n-t+1} \binom{n-v-1}{t-2}  = 1. 
\end{aligned}\end{equation}
Combining with the previous result for ${c_{n,1}^1}$ we have $c_{n,t}^1 = \mathbbm{1} \{ n \geq 2 \}$ for $t\geq 1$.

To compute ${b_{n,t}^1}$ for ${t \geq 2}$, we rely on Eq.~\eqref{eq:r_outline5} and obtain
\begin{equation}\begin{aligned}
    b_{n,t}^1 & = \frac{1}{t-1} \binom{n-1}{t-1} ^{-1} \times
    \sum_{v=1}^{n-t+1} \binom{n-v-1}{t-2} \sum_{t^\prime=1}^{t-1} c_{n-v,t'}^0 \times c_{v,1}^1 \\
    & = \binom{n-1}{t-1} ^{-1} \times \sum_{v=2}^{n-t+1} \binom{n-v-1}{t-2}  = \frac{n-t}{n-1}\,.
\end{aligned}\end{equation}

The base cases are~${b_{n,1}^1 = b_{n,2}^1 = \frac{n-2}{n-1}}$ for ${n \geq 2}$ and 
${b_{1,1}^1 = 0}$. Combining the above results yields:
\begin{align}
    b_{n,t}^1 = 
    \begin{cases}
        0                     & t = n = 1, \\
        \frac{n-2}{n-1} & t = 1, n \geq 2, \\
        \frac{n-t}{n-1} & t \geq 2,
    \end{cases}
\end{align}
from which
\begin{equation}
\label{r1_n}
    \bar{r}_1(n)  = \frac{1}{n} \sum_{t = 1}^n b_{n,t}^1 = \frac{n-2}{n(n-1)} + \frac{1}{n} \sum_{t = 2}^n \frac{n-t}{n-1} = \frac{1}{2} - \frac{1}{n(n-1)}\,.
\end{equation}
Recalling $\bar{r}_0(n)\equiv 1$ we deduce
\begin{equation}
    r_0(n)  = \bar{r}_0(n) - \bar{r}_1(n) = \frac{1}{2} + \frac{1}{n(n-1)}\,.
\end{equation}

We will follow the same procedure to compute ${\bar{r}_2(n)}$. Starting with ${c_{n,1}^k}$, we have 
\begin{equation}\begin{aligned}
    c_{n,1}^2 & = \frac{1}{n-1} c_{n-1,1}^1 + \frac{1}{n-1} \sum_{v=1}^{n-2} c_{n-v,1}^2 \cdot c_{v,1}^1 \\
    & = \frac{1}{n-1} \times \mathbbm{1} \{ n-1 > 1 \} + \frac{1}{n-1} \sum_{v=1}^{n-2} c_{n-v,1}^2 \cdot \mathbbm{1} \{ v > 1 \}\\
    & = \frac{1}{n-1} \times \mathbbm{1} \{ n > 2 \} + \frac{1}{n-1} \sum_{v=2}^{n-2} c_{v,1}^2
\end{aligned}\end{equation}
for  ${n \geq 2}$, from which 
\begin{equation}
\label{cn1:2}
 (n-1) c_{n,1}^2 = \mathbbm{1} \{ n > 2 \} + \sum_{v=2}^{n-2} c_{v,1}^2.
\end{equation}
Specializing to ${n=2}$ and ${n=3}$, we obtain ${c_{2,1}^2 = 0}$ and ${c_{3,1}^2 = 1/2}$. Massaging \eqref{cn1:2} we find
\begin{equation*}
(n-1)c_{n,1}^2 - (n-2)c_{n-1,1}^2 = c_{n-2,1}^2
\end{equation*}
for $n\geq 4$. Equivalently
\begin{equation}
c_{n,1}^2 - c_{n-1,1}^2 = - \frac{1}{n-1} (c_{n-1,1}^2 - c_{n-2,1}^2)
\end{equation}
which we iterate to find $c_{n,1}^2 - c_{n-1,1}^2 = \frac{(-1)^{n-1}}{(n-1)!}$, leading to 
\begin{align}
     c_{n,1}^2 = \sum_{j=2}^{n-1} \frac{(-1)^j}{j!}
\end{align}
for $n\geq 3$; the initial values are $c_{1,1}^2 = c_{2,1}^2 = 0$.  

Now we move to ${b_{n,t}^2}$ (since we will not be using ${c_{n,t}^2}$ for~${t\geq2}$). The evaluations for ${t=2}$ and ${t \geq 3}$ are slightly different. Using Eq.~\eqref{eq:r_outline5}, we obtain 
\begin{equation}\begin{aligned}
    b_{n,2}^2 & =  \frac{1}{n-1} \sum_{v=3}^{n-1} \mathbbm{1} \{ n-v \geq 2 \} \times \sum_{j=2}^{v-1} \frac{(-1)^j}{j!} \\
    & = \frac{1}{n-1} \sum_{j=2}^{n-3} \frac{(-1)^j}{j!} (n-j-2)\,. 
\end{aligned}\end{equation}
When ${t \geq 3}$, we get 
\begin{equation}\begin{aligned}
    b_{n,t}^2 & = \frac{1}{t-1} \binom{n-1}{t-1} ^{-1} \times
    \sum_{v=1}^{n-t+1} \binom{n-v-1}{t-2} \sum_{t^\prime=1}^{t-1} c_{n-v,t'}^1 \cdot c_{v,1}^2 \\
    & = \frac{1}{t-1} \binom{n-1}{t-1} ^{-1} \times \sum_{v=3}^{n-t+1} \binom{n-v-1}{t-2} \sum_{t^\prime=1}^{t-1} \mathbbm{1} \{ n-v \geq 2 \} \times \sum_{j=2}^{v-1} \frac{(-1)^j}{j!} \\
    & = \binom{n-1}{t-1} ^{-1} \times \sum_{j=2}^{n-t} \frac{(-1)^j}{j!} \binom{n-j-1}{t-1}.
\end{aligned}\end{equation}

We also note that ${b_{n,1}^2 = b_{n,2}^2}$, assuming ${n \geq 2}$. Combining the results for all ${t}$ yields:
\begin{align}
    & b_{n,t}^2 = 
    \begin{cases}
        \frac{1}{n-1} \sum_{j=2}^{n-3} \frac{(-1)^j}{j!} (n-j-2)      & t < 3, \\
        \binom{n-1}{t-1} ^{-1} \times \sum_{j=2}^{n-t} \frac{(-1)^j}{j!} \binom{n-j-1}{t-1} & t \geq 3.
    \end{cases}
\end{align}

Assuming ${n \geq 2}$, we compute
\begin{equation}\begin{aligned}
\label{r2_n}
    \bar{r}_2(n) & = \frac{1}{n} \sum_{t = 1}^n b_{n,t}^2 \\
    & = \frac{2}{n(n-1)} \sum_{j=2}^{n-2} \frac{(-1)^j}{j!} (n-j-2) + 
    \frac{1}{n} \sum_{t = 3}^n \binom{n-1}{t-1} ^{-1} \times \sum_{j=2}^{n-t} \frac{(-1)^j}{j!} \binom{n-j-1}{t-1} \\
   & = \frac{1}{n(n-1)} \sum_{j=2}^{n-3} \frac{(-1)^j}{(j+1)!} [-2(j+1) - j(j+1) + (n-1)n].
\end{aligned}\end{equation}

Combining \eqref{r1_n} and \eqref{r2_n} we find
\begin{equation}\begin{aligned}
    r_1(n) & = \bar{r}_1(n) - \bar{r}_2(n) \\
    & =  \frac{1}{2} - \frac{1}{n(n-1)} - \frac{1}{n(n-1)} \sum_{j=2}^{n-3} \frac{(-1)^j}{(j+1)!} [-2(j+1) - j(j+1) + (n-1)n] \\
    & = \frac{1}{2} + \sum_{j=3}^{n-2} \frac{(-1)^j}{j!} - \frac{1}{n(n-1)} \left( 1 - \sum_{j=2}^{n-3} \frac{(-1)^j}{j!} (j+2) \right)\\
    & =  \sum_{j=0}^{n-2} \frac{(-1)^j}{j!}  + \frac{1}{n(n-1)} \sum_{j=0}^{n-3} \frac{(-1)^j}{j!} (j+2)
\end{aligned}\end{equation}
as stated in Theorem \ref{the:finite_n}. 

To compute $\ell_0(n)$ and~$\ell_1(n)$, we use the quantities $c_{n,t}^k$. We begin with recurrence
\begin{align}
     c_{n,t}^k = \sum_{v=1}^{n-t+1} h_{n,t}(v) ~\mathbbm{E}_{t'<t} \left[ 1 - (1 - c_{n-v,t'}^{k-1})(1 - c_{v,1}^k) \right],
\end{align}
which we specialize to ${t=2}$ and use ${c_{n,1}^k = c_{n,2}^k}$ to get 
\begin{align}
\label{cn1:k}
     c_{n,1}^k = \frac{1}{n-1} \sum_{v=1}^{n-1} \left[ 1 - (1 - c_{n-v,1}^{k-1})(1 - c_{v,1}^k) \right].
\end{align}

We have ${c_{n,1}^0=1}$. Specializing \eqref{cn1:k} to ${k=1}$ and keeping in mind that ${c_{1,1}^1 = 0}$, we get 
\begin{align}
     c_{n,1}^1 = \frac{1}{n-1} \sum_{v=1}^{n-1} \left[ 1 - (1 - 1)(1 - c_{v,1}^1) \right] = 1
\end{align}
for ${n \geq 2}$, i.e., ${c_{n,1}^1 = 1 - \delta_{n,1}}$. Specializing \eqref{cn1:k} to ${k=2}$, keeping in mind that ${c_{1,1}^2 = 0}$, and using previous results we find
\begin{equation}
c_{n,1}^2 = \frac{1}{n-1} \sum_{v=1}^{n-1} \left[ 1 - \delta_{n-v,1}(1 - c_{v,1}^2) \right]
\end{equation}
for~${n \geq 2}$. Equivalently
\begin{equation}
c_{n,1}^2 = \frac{n-2}{n-1} + \frac{1}{n-1} c_{n-1,1}^2,
\end{equation}
which is solved subject to ${c_{1,1}^2 = 0}$ to yield~${c_{n,1}^2 = 1 - \frac{1}{(n-1)!}}$.

We now consider ${c_{n,t}^k}$ with ${t \geq 2}$. We have ${c_{n,t}^0 = 1}$ and ${c_{1,t}^1 = 0}$, while for~${n \geq 2}$
\begin{align}
     c_{n,t}^1 = \sum_{v=1}^{n-t+1} h_{n,t}(v) ~\mathbbm{E}_{t'<t} \left[ 1 - (1 - 1)(1 - c_{v,1}^k) \right] = 1.
\end{align}
Thus,~${c_{n,t}^1 = 1 - \delta_{n,1}}$. For~${k=2}$,
\begin{equation}\begin{aligned}
    c_{n,t}^2 & = \sum_{v=1}^{n-t+1} h_{n,t}(v) ~\mathbbm{E}_{t'<t} \left[ 1 - (1 - c_{n-v,t'}^{k-1}) \frac{1}{(v-1)!} \right] \\
    & = \delta_{t,2} \left[ 1 - \frac{1}{(n-1)!} \right] + (1-\delta_{t,2}) \\
    & = 1 - \frac{\delta_{t,2}}{(n-1)!}.
\end{aligned}\end{equation}

To determine ${b_{n,t}^k}$, we rely on the recurrence
\begin{align}
\label{bnt}
     b_{n,t}^k = \sum_{v=1}^{n-t+1} h_{n,t}(v) ~\mathbbm{E}_{t'<t} \left[ c_{n-v,t'}^{k-1} \cdot c_{v,1}^k + \left( 1 - c_{n-v,t'}^{k-1} \right) \cdot b_{v,1}^k \right]
\end{align}
valid when $t\geq 2$. The case of ${t=1}$ is deduced from ${b_{n,1}^k = b_{n,2}^k}$. Equation \eqref{bnt} simplifies to 
\begin{equation}\begin{aligned}
    b_{n,t}^1 & = \sum_{v=1}^{n-t+1} h_{n,t}(v) ~\mathbbm{E}_{t'<t} \left[ 1 \cdot (1 - \delta_{v,1}) + \left( 1 - 1 \right) \cdot b_{v,1}^1 \right] \\
    & = 1 - h_{n,t}(1) = \frac{n-t}{n-1} 
\end{aligned}\end{equation}
when $k=1$. 
For~${k=2}$, we have
\begin{align}
 \label{bnt:2}
    b_{n,t}^2 = \sum_{v=1}^{n-t+1} h_{n,t}(v) \left[ (1 - \delta_{n-v,1}) \left( 1 - \frac{1}{(v-1)!} \right) + \delta_{n-v,1} \cdot b_{v,1}^2 \right].
\end{align}
Specializing \eqref{bnt:2} to ${t=2}$ and using~${b_{n,1}^k = b_{n,2}^k}$ we obtain
\begin{align}
     b_{n,1}^2 = \frac{1}{n-1} b_{n-1,1}^2 + \frac{n-2}{n-1} - \frac{1}{n-1} \sum_{v=1}^{n-2} \frac{1}{(v-1)!}\,.
\end{align}
The solution to this recurrence reads
\begin{align}
     b_{n,1}^2 = \frac{1}{(n-1)!} \sum_{j=1}^{n-2} \left[ j! \left( j - \sum_{v=1}^j \frac{1}{(v-1)!} \right) \right].
\end{align}

For~${t \geq 3}$, we have
\begin{equation}\begin{aligned}
    b_{n,t}^2 & = \sum_{v=1}^{n-t+1} h_{n,t}(v) \left[ 1 - \frac{1}{(v-1)!} \right] \\
    & = 1 - \binom{n-1}{t-1}^{-1} \sum_{v=1}^{n-t+1} \binom{n-v-1}{t-2} \frac{1}{(v-1)!}.
\end{aligned}\end{equation}

We can now compute ${\bar{\ell}_k(n)}$ for ${k=1}$ and $k=2$. In the first case we get 
\begin{equation}\begin{aligned}
    \bar{\ell}_1(n) & = \frac{1}{n} \sum_{t=1}^n b_{n,t}^1 \\
    & = \frac{1}{n} \left[ \frac{n-2}{n-1} + \sum_{t=2}^n \frac{n-t}{n-1} \right] = \frac{1}{2} - \frac{1}{n(n-1)}
\end{aligned}\end{equation}
coinciding, as it should, with $\bar{r}_1(n)$ given by  \eqref{r1_n}, and thus providing the consistency check. 

When~${k=2}$, we have
\begin{equation}\begin{aligned}
    \bar{\ell}_2(n) & = \frac{1}{n} \sum_{t=1}^n b_{n,t}^2 \\
    & = \frac{2}{n!} \sum_{j=1}^{n-2} \left[ j! \left( j - \sum_{v=1}^j \frac{1}{(v-1)!} \right) \right] + \frac{1}{n} \sum_{t=3}^n \left[ 1 - \binom{n-1}{t-1}^{-1} \sum_{v=1}^{n-t+1} \binom{n-v-1}{t-2} \frac{1}{(v-1)!} \right].
\end{aligned}\end{equation}
Thus
\begin{equation}\label{eq:th.7.proof.l1}
\begin{aligned}
    \ell_1(n) = & \,\, \bar{\ell}_1(n) - \bar{\ell}_2(n) \\
    = & -\frac{1}{2} + \frac{3}{n} - \frac{1}{n-1} - \frac{2}{n!} \sum_{j=1}^{n-2} \left[ j! \left( j - \sum_{v=1}^j \frac{1}{(v-1)!} \right) \right] \\ & + \frac{1}{n} \sum_{t=3}^n \binom{n-1}{t-1}^{-1} \sum_{v=1}^{n-t+1} \binom{n-v-1}{t-2} \frac{1}{(v-1)!}.
\end{aligned}
\end{equation}

To simplify this expression, we note that
\begin{equation}
    \frac{1}{n} \sum_{t=3}^n \binom{n-1}{t-1}^{-1} \sum_{v=1}^{n-t+1} \binom{n-v-1}{t-2} \frac{1}{(v-1)!}
    = \sum_{t=3}^n \sum_{v=1}^{n-t+1} \frac{(n-v-1)!}{n!} (t-1) \binom{n-t}{v-1},
\end{equation}
the right hand side of which with ${\delta=n-t}$ simplifies as follows:
\begin{align*}
    & \sum_{\delta=0}^{n-3} \sum_{v=1}^{\delta+1} \frac{(n-v-1)!}{n!} (n-\delta-1) \binom{\delta}{v-1} \\
    & = \sum_{\delta=0}^{n-3} \sum_{v=1}^{\delta+1} \sum_{c=1}^{n-\delta-1} \frac{(n-v-1)!}{n!} \binom{\delta}{v-1} \\
    & = \sum_{v=1}^{n-2} \frac{(n-v-1)!}{n!} \sum_{\delta=v-1}^{n-3} \sum_{c=1}^{n-\delta-1} \binom{\delta}{v-1} \\
    & = \sum_{v=1}^{n-2} \frac{(n-v-1)!}{n!} \left( \sum_{c=2}^{n-v} \sum_{\delta=v-1}^{n-c-1} \binom{\delta}{v-1} + \sum_{\delta=v-1}^{n-3} \binom{\delta}{v-1} \right) \\
    & = \sum_{v=1}^{n-2} \frac{(n-v-1)!}{n!} \left( \sum_{c=2}^{n-v} \binom{n-c}{v} + \binom{n-2}{v} \right) \\
    & = \sum_{v=1}^{n-2} \frac{(n-v-1)!}{n!} \left( \binom{n-1}{v+1} + \binom{n-2}{v} \right) \\
    & = \sum_{v=2}^{n-1} \frac{1}{v!} - \frac{1}{n(n-1)} \sum_{v=0}^{n-3} \frac{1}{v!}.
\end{align*}
Plugging this into~\eqref{eq:th.7.proof.l1} and simplifying further leads to
\begin{align}
    \ell_1(n) =& -\frac{1}{2} + \frac{1}{n} - \frac{1}{n-1}
    + \frac{2}{n!} + \sum_{v=2}^{n-1} \frac{1}{v!} - \frac{1}{n(n-1)} \sum_{v=0}^{n-3} \frac{1}{v!} + \frac{2}{n!} \sum_{j=1}^{n-2} \sum_{v=0}^{j-1} \frac{j!}{v!},
\end{align}
the extraction of leading terms from which yields
\begin{align}
    \ell_1(n) &= e - \frac{5}{2} - \frac{e+1}{n(n-1)}
    + \frac{2}{n!} \sum_{v=0}^{n-3} \sum_{j=v+1}^{n-2} \frac{j!}{v!} + O \left( \frac{1}{n!} \right) \nonumber\\
    &= e - \frac{5}{2} + \frac{e-1}{n^2} + O \left( \frac{1}{n^3} \right).
\end{align}

\appendix

\section{Traveling wave analysis for Conjecture~\ref{con:layer_dist_tail}}
\label{ssec:large_k_proof_outline}

In this section, we describe the traveling wave method that we use to obtain the large-${k}$ tail of the layer distribution~\eqref{eq:layer_tail} stated in Conjecture~\ref{con:layer_dist_tail}.

Traveling wave phenomena underlying the propagation of the front from a stable state into an unstable state were originally found in the context of partial differential equations~\cite{Bramson,van2}, and later also in discrete settings involving recurrences or difference-differential equations~\cite{Krapivsky2000, Majumdar2002, Krapivsky2012}. One often finds a class of solutions propagating at velocities below or above a certain threshold, and argues that, for most natural (usually sufficiently steep) initial conditions, the unique solution that advances at the extremal threshold velocity is selected. More precisely, the front propagates ballistically in the leading order, with a sub-leading logarithmic correction~\cite{Bramson, Brunet_1997, Derrida_2023}. This extremal selection mechanism has been rigorously proven for a limited class of equations~\cite{Bramson, Derrida_2023, jonckheere2026waves}.

Here, we first show that a traveling wave indeed emerges from the proven recurrence in Eq.~\eqref{eq:pk2}. We then apply the velocity selection principle to the resulting otherwise-undetermined dispersion relation. We also conjecture the associated logarithmic correction, and a common asymptotic scaling of~$\bar\ell_k$ and~$\ell_k$. We support each conjectural step by numeric evidence. Our approach follows closely the application of the traveling wave method in~\cite{Krapivsky2012} to the derivation of the asymptotic height of a cascade tree.

\begin{figure}
    \centering
    \includegraphics[width=0.95\textwidth,trim = 10 0 0 0]{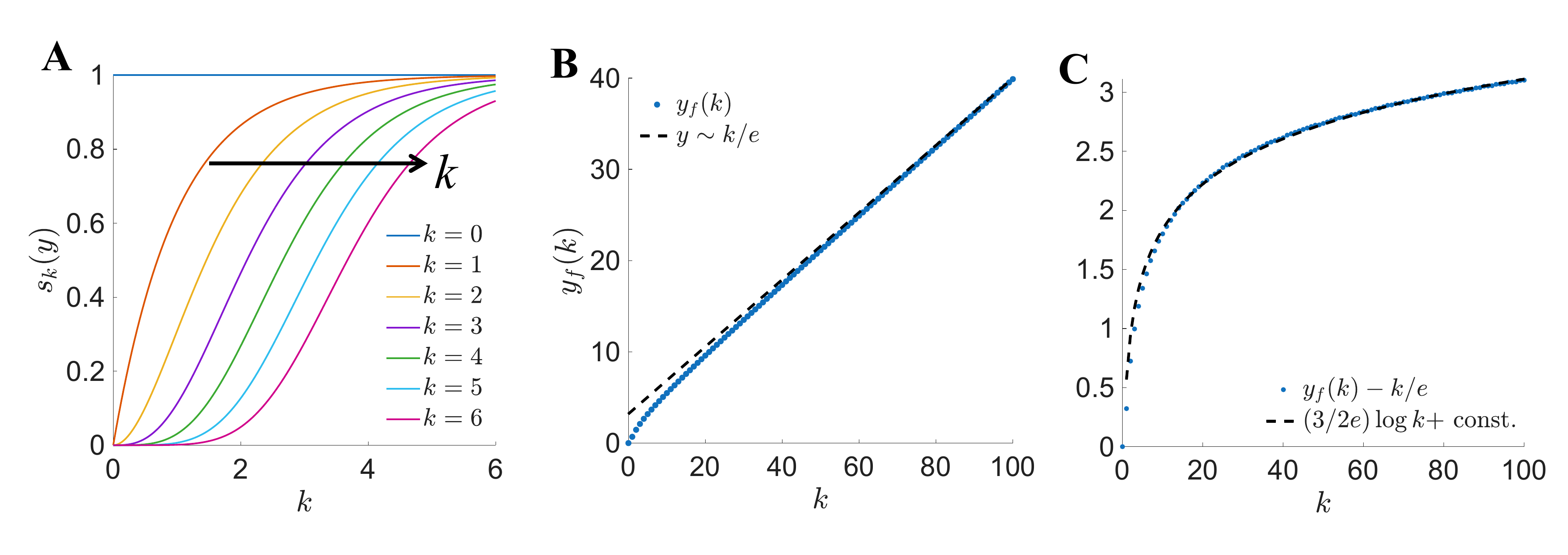}
    \caption{\textbf{(A)}~The auxiliary functions~${s_k(y)}$ converge to a traveling wave at~${k \rightarrow \infty}$. The arrow indicates the direction of growing~${k}$. \textbf{(B)}~For sufficiently large~${k}$, the wavefront~${y_f(k)}$ grows with~$k$ as~${y_f \sim k/e}$. The wavefront position is computed as $y_f(k)=\mathrm{argmin}_{y\ge 0}~|s_k(y)-1/2| $. \textbf{(C)}~The sublinear correction to wavefront for $\bar\ell_k$.
    }
    \label{FF}
\end{figure}

We first recall that Theorem~\ref{the:pk} provides the cumulative layer distribution
\begin{equation}
\label{eq:bar_ell_k_appendix}
\bar{\ell}_k=\int_0^\infty s_k(y)e^{-y}dy
\end{equation}
in terms of auxiliary functions $\{s_k(y)\}_{k\ge 0}$, defined recursively by
\begin{equation}
    \label{eq:pk2_appendix}
    s_0(y) = 1,
    \quad s_k(y) = 1 - \exp\!\left[ -\int_0^y s_{k-1}(y^\prime) \, dy^\prime \right] \quad \forall k \geq 1.
\end{equation}
Let us extend the domain of~${s_k(y)}$ to the entire real line by defining~${s_k(y) = 0}$ for~${y < 0}$. With this extension, the recursion \eqref{eq:pk2_appendix} becomes
\begin{align}
    \label{E_gk_recursion_modified}
   {s_0(y) = \Theta(y)}, \quad s_k(y) = 1 - \exp \left[ -\int_{-\infty}^y dy'~ s_{k-1}(y') \right],
\end{align}
with $\Theta(y):=\mathds{1}\{y\ge 0\}$ denoting the Heaviside step function. Differentiating, we obtain an equivalent differential form of \eqref{E_gk_recursion_modified} for $k\ge 1$:
\begin{align}
    \label{E_gk_recursion_modified2}
    & \frac{ds_k(y)}{dy} = s_{k-1}(y) \bigl( 1 - s_k(y) \bigr).
\end{align}
We now introduce a traveling wave ansatz with a fixed waveform~${s(y)}$, and a wavefront position~${y_f(k)}$ for large~$k$:
\begin{align}
    \label{E_traveling_wave_ansatz}
    & \lim_{k \rightarrow \infty} \bigl[ s_k(y) - s\bigl(y-y_f(k)\bigr) \bigr] = 0.
\end{align}
The numerical evaluation of~${s_k(y)}$ supports this ansatz convergence in Fig.~\ref{FF}(A). We further assume a linear growth in wavefront position to leading order,~${y_f(k) \sim vk}$, with an asymptotic wave velocity~${v>0}$. This linear growth to leading order is confirmed by the numerical evaluation in Fig.~\ref{FF}(B). 

Under these numerically validated assumptions, substituting~${s_k(y)\approx s\bigl(y-y_f(k)\bigr)}$ into~\eqref{E_gk_recursion_modified2} and taking the large~${k}$ limit yields
\begin{align}
\label{eq:dsdy_wave}
    & \frac{ds(y)}{dy} = s(y+v)(1-s(y)).
\end{align}
Since for every~${k}$ we have~${s_k(y<0) = 0}$, the asymptotic waveform~${s(y)}$ must satisfy~${s(y)\rightarrow0}$ as~${y \rightarrow -\infty}$. Thus, for sufficiently large negative~${y}$ we may approximate~${1-s(y)\approx 1}$ in~\eqref{eq:dsdy_wave}, leading to the following delay differential equation
\begin{align}
    & \frac{ds(y)}{dy} = s(y+v),
\end{align}
which admits solutions of the form
\begin{equation}
  s(y) \sim e^{ay},
\end{equation}
where $a$ is a solution of the dispersion relation
\begin{equation}
  e^{av} = a.
\end{equation}
Solving the latter for~${v}$ gives
\begin{equation}
v(a) = a^{-1} \log a,
\end{equation}
and hence an upper bound on wave velocity
\begin{equation}
  v \leq 1/e. 
\end{equation}

We next assert that the actual wave velocity attains its extremal value~${v^* = 1/e}$. This unusual feature of the extremal velocity being ``selected'' has been observed in a large class of similar problems involving nonlinear traveling wave solutions with Heaviside initial conditions~\cite{Majumdar2002}. Its rigorous characterization remains an open problem. We confirm the value $1/e$ in the numerical evaluation as displayed in Fig.~\ref{FF}(B).

Next, we note that in a similar large class of recursions~\cite{Krapivsky2000, Majumdar2002, Krapivsky2012}, the subleading correction to the wavefront position has the universal form~${\frac{3}{2a^*} \log k}+O(1)$, where~$a^*$ solves the dispersion relation at the extremal velocity value~$v^*$. In our case,~${a^*=e}$, corresponding to $e^{av^*}=a$ with~${v^*=1/e}$, so we conjecture 
\begin{align}
\label{E_front}
    & y_f(k) = \frac{k}{e} + \frac{3}{2e} \log k + O(1),
\end{align}
with the sublinear correction supported by the numerics in Fig.~\ref{FF}(C).

Substituting the last expression into Eq.~\eqref{eq:bar_ell_k_appendix} for~${\bar{\ell}_k}$ (Theorem~\ref{the:pk}), we obtain
\begin{equation}
    \bar{\ell}_k = e^{-y_f(k)} \int_{-y_f(k)}^\infty s_k\left(y'+y_f(k)\right) \, e^{-y'} dy'.
\end{equation}
At~${k\rightarrow\infty}$, the integral converges to a constant ${\int_{-\infty}^\infty s(y') \, e^{-y'} dy'}$. Thus~${\bar{\ell}_k \sim e^{-y_f(k)}}$, and by Eq.~\eqref{E_front} we conclude that
\begin{align}
\label{E_traveling_wave_1}
    \bar{\ell}_k & = \exp \left( -\frac{k}{e} - \frac{3}{2e} \log k - O(1) \right).
\end{align}

To complete the argument in favor of Conjecture~\ref{con:layer_dist_tail}, we further assert that~${\ell_k \sim \bar{\ell}_k}$, so that~${\ell_k}$ has the asymptotic behavior given in Eq.~\eqref{E_traveling_wave_1}. Although one can not infer this result directly from~\eqref{E_traveling_wave_1} due to the unknown~${O(1)}$ term, it is supported by the numerical evaluation of the~${\ell_k}$ at large~${k}$ in Fig.~\ref{fig:RLdist_numerical}(D).

\bibliography{bib}

\end{document}